\documentclass[11pt]{article}

\usepackage{amsmath}
\usepackage{mathtools}
\usepackage{amssymb}
\usepackage{latexsym}
\usepackage{enumitem}
\usepackage{mathrsfs}
\usepackage{comment}
\usepackage{pifont}
\usepackage{enumitem}
\usepackage{bbold}
\usepackage{graphicx}
\usepackage{color}
\usepackage{hyperref}
\hypersetup{
    colorlinks=true,
    linkcolor=blue,
    filecolor=magenta,      
    urlcolor=cyan,
    citecolor=magenta
}

\newtheorem{assumption}{Assumption}

\def\qed{ \ \vrule width.2cm height.2cm depth0cm\smallskip}
\newenvironment{proof}{\noindent {\bf Proof.\/}}{$\qed$\vskip 0.1in}

\newcommand{\ol}{\overline}
\newcommand{\ul}{\underline}

\newcommand{\ba}{\begin{array}}
\newcommand{\ea}{\end{array}}
\newcommand{\be}{\begin{equation}}
\newcommand{\ee}{\end{equation}}
\newcommand{\bea}{\begin{eqnarray}}
\newcommand{\eea}{\end{eqnarray}}
\newcommand{\beaa}{\begin{eqnarray*}}
\newcommand{\eeaa}{\end{eqnarray*}}

\def\dbE{\mathbb{E}}
\def\dbF{\mathbb{F}}

\def\dbL{\mathbb{L}}

\def\dbP{\mathbb{P}}
\def\dbR{\mathbb{R}}
\def\dbS{\mathbb{S}}

\def\dbQ{\mathbb{Q}}

\def\a{\alpha}
\def\b{\beta}
\def\g{\gamma}
\def\d{\delta}
\def\e{\varepsilon}

\def\l{\lambda}
\def\m{\mu}
\def\n{\nu}
\def\si{\sigma}
\def\t{\tau}
\def\f{\varphi}
\def\th{\theta}
\def\o{\omega}

\def\G{\Gamma}

\def\L{\Lambda}

\def\O{\Omega}
\def\cA{{\cal A}}
\def\cB{{\cal B}}
\def\cC{{\cal C}}

\def\cE{{\cal E}}
\def\cF{{\cal F}}

\def\cN{{\cal N}}

\def\cP{{\cal P}}

\def\cS{{\cal S}}
\def\cT{{\cal T}}

\def\cW{{\cal W}}

\def\no{\noindent}

\def\q{\quad}
\def\qq{\qquad}

\def\pa{\partial}
\def\cd{\cdot}

\def\tr{\hbox{\rm tr}}

\def\qed{ \hfill \vrule width.25cm height.25cm depth0cm\smallskip}

\newcommand{\basa}{\begin{assumption}}
\newcommand{\easa}{\end{assumption}}

\newcommand{\bas}{\begin{assum}}
\newcommand{\eas}{\end{assum}}

\def\pa{\partial}

 \def\cd{\cdot}

\def\supp{\hbox{\rm supp$\,$}}

\def\tr{\hbox{\rm tr$\,$}}

\def\bx{{\bf x}}
\def\bX{{\bf X}}
\def\bY{{\bf Y}}
\def\bA{{\bf A}}
\def\bM{{\bf M}}
\def\bz{{\bf z}}
\def\be{{\bf e}}
\def\bo{{\boldsymbol \o}}
\def\balpha{{\boldsymbol \a}}
\def\bsigma{{\boldsymbol \si}}
\def\bgamma{{\boldsymbol \g}}
\def\ba{{\bf a}}

\def\1{{\bf 1}}
\def\by{{\bf y}}

\def\:{\!:\!}
\def\reff#1{{\rm(\ref{#1})}}
\def \proof{{\noindent \bf Proof\quad}}

\def \bQ{{\bf Q}}

at 9pt

\begin{document}

\newtheorem{thm}{Theorem}[section]
\newtheorem{lem}[thm]{Lemma}
\newtheorem{cor}[thm]{Corollary}
\newtheorem{prop}[thm]{Proposition}
\newtheorem{rem}[thm]{Remark}
\newtheorem{eg}[thm]{Example}
\newtheorem{defn}[thm]{Definition}
\newtheorem{assum}[thm]{Assumption}

\numberwithin{equation}{section}

\title{\bf{A finite-dimensional approximation for partial differential equations on Wasserstein space}\footnote{The author is grateful for the financial support from the Chaires FiME-FDD and Financial Risks of the Louis Bachelier Institute. He also thanks Nizar Touzi and Fabrice Mao Djete for their insightful comments.}}
\author{Mehdi Talbi\footnote{Laboratoire de Probabilités, Statistique et Modélisation, Université Paris-Cité, France, talbi@lpsm.paris}}
\date{\today}

\maketitle

\begin{abstract}
This paper presents a finite-dimensional approximation for a class of partial differential equations on the space of probability measures. These equations are satisfied in the sense of viscosity solutions. The main result states the convergence of the viscosity solutions of the finite-dimensional PDE to the viscosity solutions of the PDE on Wasserstein space, provided that uniqueness holds for the latter, and heavily relies on an adaptation of the Barles \& Souganidis monotone scheme \cite{BarSou} to our context, as well as on a key precompactness result for semimartingale measures. We illustrate our convergence result with the example of the Hamilton-Jacobi-Bellman and Bellman-Isaacs equations arising in stochastic control and differential games, and propose an extension to the case of path-dependent PDEs.
\end{abstract}

\no{\bf MSC2010.} 35D40, 35R15, 60H30, 49L25.

\vspace{3mm}
\no{\bf Keywords.}  Viscosity solutions, Wasserstein space, path-dependent PDEs, mean field control.

\section{Introduction}

Since Lasry \& Lions \cite{LL} and Caines, Huang \& Malhamé \cite{CHM} introduced mean field games, partial differential equations on the space of probability measures have become a popular tool to study systems of interacting agents. In \cite{CDLL}, Cardaliaguet, Delarue, Lasrly \& Lions used the master equation to prove the convergence of the $N$-players game to the corresponding mean field problem, provided that the value function of the problem is smooth. Since then, an extensive part of the literature has been dedicated to these convergence issues: see e.g. Cardaliaguet \cite{Cardaliaguet}, Bayraktar, Cecchin, Cohen \& Delarue \cite{BCCD, BCCD2}, Cecchin \& Pelino \cite{CecPel}, Cecchin, Dai Pra, Fischer \& Pelino \cite{CDFP}, Laurière \& Tangpi \cite{LauTan}, Djete \cite{Djete} or Doncel, Gast \& Gaujal \cite{DGG}.

In the context of mean field control, an extensive literature focuses on the convergence issue, see e.g.\ Bayraktar, Cecchin \& Chakraborty \cite{BayCha}, Cardaliaguet, Daudin, Jackson \& Souganidis \cite{CDJS}, Cardaliaguet \& Souganidis \cite{CarSou},  Cavagnari, Lisini, Orrieri \& Savare \cite{CLOS}, Cecchin \cite{Cecchin},  Djete, Possamai \& Tan  \cite{DPT},  Fischer \& Livieri \cite{FisLiv},  Fornasier, Lisini, Orrieri \& Savare \cite{FLOS}, Germain, Pham \& Warin \cite{GPW}, and Lacker \cite{LackerPropagation}. In particular, by utilizing some strong regularity of the value function, \cite{GPW, CDJS} obtained certain rate of convergence. 

The contributions the most relevant to our work are the ones of Gangbo, Mayorga \& Swiech \cite{GMS}, Mayorga \& Swiech \cite{mayorga2023finite} and Talbi, Touzi \& Zhang \cite{TTZ3}. In the first two papers, the authors develop a finite-dimensional approximation for Hamilton-Jacobi-Bellman equations with uncontrolled volatility, in the sense of viscosity solutions (defined via lifting on the Hilbert space of square integrable random variables). The second paper introduces an approximation of the obstacle problem on Wasserstein space, which characterizes the mean field optimal stopping problem (see \cite{TTZ, TTZ2}).  

Our objective is to find a finite-dimensional approximation for a general class of PDEs on Wasserstein space, satisfied in the sense of viscosity. We use the notion of viscosity solutions developed by Wu \& Zhang \cite{WZ}, which is intrinsic and allows for path-dependent PDEs (i.e., the solutions of the equations depend on the time and on a probability measure on the space of continuous paths of $\dbR^d$). This class of equations covers in particular equations arising in mean field stochastic control (including the case of controlled volatility). We intend to find a finite-dimensional PDE whose viscosity solutions converge to viscosity solutions of the equation on Wasserstein space. More precisely, we follow the methodology of \cite{TTZ3} and adapt Barles \& Souganidis' monotone scheme \cite{BarSou} to our context, by proving that the semi-relaxed limits of the viscosity super/subsolutions of the finite-dimensional PDE are viscosity super/subsolutions of the PDE on Wasserstein space as the dimension goes to infinity. If uniqueness holds for the latter equation, this implies the convergence of viscosity solutions. 

We propose an application of this result to the convergence of the value functions of finite dimensional control problem to the value function of mean field control problems. In particular, we show that this convergence holds under quite weak regularity assumption on the coefficients and rewards of the problem. 

An important feature that ensures the convergence relies on the choice of the set of the test functions for the finite dimensional problem. Similarly to the viscosity theory developed for path-dependent PDEs (see Ekren, Keller, Ren, Touzi \& Zhang \cite{EKTZ, ETZ1, ETZ2, RTZ}; see also Guo, Zhang \& Zhuo \cite{zhang2014monotone} and Ren \& Tan \cite{ren2017convergence}), we only require test functions to be tangent to the super/subsolution through the mean, whereas the tangency is pointwise in the standard literature. This tangency in expectation can be seen as the finite-dimensional counterpart of our set of test functions on Wasserstein space. We emphasize that this choice of test functions is crucial even in the Markovian case: indeed, by considering only points along the trajectories of the state process, we are able to apply a key propagation of chaos-like result and to derive essential estimates to prove the main convergence theorem.

One of the main advantages of our methodology is that it is able to handle very similarly the Markovian and path-dependent cases. Although we discuss in more detail the Markovian equations for the sake of clarity, we insist on the fact that the proofs for non-Markovian equations are (almost) the same, and that most change only lie in some definitions and notations.  We also emphasize that our results apply to a large class of equations, exceeding the scope Hamilton-Jacobi-Bellman equations on Wasserstein space and including for examples equations with non-convex Hamiltonians (such as Bellman-Isaacs equations). 


The paper is organized as follows. In Section \ref{sect-equation}, we present the class of equations on Wasserstein space and the notion of viscosity solutions. In Section \ref{sect-approximation}, we introduce the finite-dimensional approximation and state the main result of the paper, the convergence theorem, which we apply in Section \ref{sect-examples} to mean field stochastic control. Section \ref{sect-PPDE} extends our results to the case of path-dependent PDEs. Sections \ref{sect-proof} and \ref{sect-propagation} are respectively dedicated to the proof of the convergence theorem and of the precompactness result. 

\no {\bf Notations.} Let $(E, \cA)$ be a measurable space endowed with a metric $d$. We denote by $\cP(E,\cA)$ the set of probability measures on $(E,\cA)$, and by $\cP_p(E,\cA)$ its subset of  $p$-integrable probability measures, $p \ge 1$, equipped with the Wasserstein distance defined by 
$$ \cW_p(\m, \n) := \inf_{\dbQ \in \Pi(\m, \n)} \Big(\int_{E \times E} d^p(x, y) \dbQ(x, y) \Big)^{1/p} \q \mbox{for all $(\m, \n) \in \cP_p(E, \cA)^2$,} $$
where $\Pi(\m, \n)$ is the set of couplings of $\m$ and $\n$. When $\cA = \cB(E)$, the Borel $\si$-algebra of $E$, we simply write $\cP(E)$ and $\cP_p(E)$. We denote by $\supp(\m)$ the support of $\m \in \cP(E, \cA)$, defined as the smallest closed set $\cC \subset E$ s.t. $\m(\cC^c) = 0$. Given a random variable $Z$ and a probability $\dbP$, we denote by $\dbP_Z:=\dbP\circ Z^{-1}$ the law of $Z$ under $\dbP$. We shall sometimes write $\langle \m, f \rangle := \int f d\m$. The space of the $d\times d$ real valued symmetric matrices is denoted by $\dbS_d$, and $\dbS_{d \times N}^{\rm D}$ denotes the set of blockwise diagonal matrices of the form ${\rm Diag}(A_1, \dots, A_N)$, where each $A_i \in \dbS_d$. For vectors $x, y\in \dbR^d$ and matrices $A, B$, denote $ x\cd y:=\sum_{i=1}^d x_iy_i$  and $A:B:= \tr(A B^\top)$. Given $\bx := (x_1, \cdots, x_N) \in E^N$, we denote by $\m^N(\bx) := \frac 1 N \sum_{i=1}^N \d_{x_i} \in \cP(E)$, and $\cP^N(E) := \{\m^N(\bx) : \bx \in E^N\}$. For $p \ge 1$, we also write
\beaa
\lVert \bx \rVert_p := \Big( \frac 1 N \sum_{i=1}^N | x_i |^p \Big)^{1/p} \ \mbox{and} \ \lVert \m \rVert_p := \Big( \int_{\dbR^d} \lvert x \rvert^p \m(dx) \Big)^{1/p} \q \mbox{for all $\bx \in \dbR^{d \times N}$ and $\m \in \cP_p(\dbR^d).$}
\eeaa
We shall also write “LSC" (resp. “USC") for “lower (resp. upper) semi-continuous".

\section{Viscosity solutions of partial differential equations on Wasserstein space}\label{sect-equation}

\subsection{Differentiability on Wasserstein space}

For $t \in [0,T)$, we denote 
$$ \bQ_t := [t,T) \times \cP_2(\dbR^d) \q \mbox{and} \q \ol \bQ_t := [t,T] \times \cP_2(\dbR^d). $$

 \begin{defn}\label{C12S} 
 Fix $t \in [0,T)$.
 
\no {\rm (i)} $u : \ol \bQ_t \longrightarrow \dbR$ has a functional linear derivative if there exists $\d_m u:  \ol \bQ_0 \times \dbR^d \to \dbR$ satisfying, for any $s\in [t, T]$ and $\m, \n \in \cP_2(\dbR^d)$,
\beaa
u(t,\n)-u(t, \m) = \displaystyle\int_0^1 \int_{\dbR^d} \d_m u(s, \l \n + (1-\l)\m, x)(\n-\m)(dx)d\l,
\eeaa
and $\d_m u$ has quadratic growth in $x \in \dbR^d$, locally uniformly in $(s,m) \in \ol \bQ_t$.

\no{\rm (ii)} We denote by $C^{1,2}_b(\ol \bQ_t)$ the set of bounded functions $u: \ol \bQ_t \to \dbR$ such that $\pa_t u$, $\d_m u$, $\pa_x \d_m u$, $\pa_{xx}^2\d_m u$ exist and are continuous and bounded in all their variables.
 \end{defn}

\subsection{Partial differential equation on Wasserstein space}

Let $F : [0,T] \times \cP_2(\dbR^d) \times \dbR \times \dbL_2^0(\dbR^d, \dbR^d) \times \dbL_2^0(\dbR^d, \dbS_d) \longrightarrow \dbR$, where $\dbL_2^0(\dbR^d,\dbR^d)$ (resp. $\dbL_2^0(\dbR^d, \dbS_d)$) denotes the set of Borel measurable functions from $\dbR^d$ to $\dbR^d$ (resp. $\dbS_d$) with quadratic growth, and $g : \cP_2(\dbR^d) \longrightarrow \dbR$.  We consider the following equation:
\bea\label{pde}
-\pa_t u(t,\m) - F\big(t, \m, u(t,\m), \pa_x \d_m u(t,\m,\cdot), \pa_{xx}^2 \d_m u(t,\m,\cdot) \big) = 0, \ u|_{t=T} = g, \ \mbox{$(t,\m) \in \bQ_0$.} 
\eea
The following assumptions on $F$ will be crucial to guarantee the existence of a finite-dimensional approximation to the solution of \eqref{pde}:

\bas\label{assum-operator}
{\rm (i)} $F$ is continuous in the following sense: 
$$F(t^n,\m^n, y^n, Z^n, \G^n) \longrightarrow F(t, \m, y, Z, \G) \ \mbox{as $(t^n, \m^n, y^n, Z^n, \G^n) \underset{n \to \infty}{\longrightarrow} (t, \m, y,  Z,  \G)$,}$$
where the convergence of $\m^n$ to $\m$ is in $(\cP_2(\dbR^d), \cW_2)$ and the convergence of $(Z^n, \G^n)$ to $(Z, \G)$ is pointwise for $Z,Z^n \in C^0(\dbR^d, \dbR^d)$ and $ \G, \G^n \in C^0(\dbR^d, \dbS_d)$ with quadratic growth uniformly in $n$. 
\\
{\rm (ii)} For all $(t,\m, y, Z, \G) \in [0,T] \times \cP_2(\dbR^d) \times \dbR \times  \dbL_2^0(\dbR^d, \dbR^d) \times \dbL_2^0(\dbR^d, \dbS_d) $, we have 
$$ F(t, \m, y, Z, \G) = F(t, \m, y, Z', \G') $$
for all $Z', \G'$ s.t. $Z'|_{\tiny \supp(\m)}, \G'|_{\tiny \supp(\m)} = Z|_{\tiny \supp(\m)}, \G|_{\tiny \supp(\m)}$.
\eas


\subsection{Viscosity solutions}

Let $L$ be a positive constant, $\O := C^0([0,T], \dbR^d)$, $X$ be the canonical process on $\O$ and $\dbF := \{ \cF_t \}_{t \in [0,T]}$ be the corresponding filtration.
\begin{defn}
We denote by $\cP_L$ the set of measures $\dbP \in \cP_2(\O, \cF_T)$ s.t.\ $X$ is a $\dbP$-semimartingale with drift and diffusion characteristics uniformly bounded by $L$. For $(t,\m) \in \bQ_0$, we also define $\cP_L(t,\m) := \{ \dbP \in \cP_L : \dbP_{X_t} = \m \}$.
\end{defn}
We define the viscosity neighborhood of $(t,\m) \in \bQ_0$ by
$$ \cN_\d (t,\m) := \{(s, \dbP_{X_s}) :  s \in [t, t+\d], \dbP \in \cP_L(t,\m) \}, $$
Note that the constant $L$ may be chosen arbitrarily large, and has for unique purpose to ensure the $\cW_2$-compactness of $\cP_L(t,\m)$ for all $(t, \m)$ and, therefore, the one of $\cN_\d(t,\m)$ (see Wu \& Zhang \cite[Lemma 4.1]{WZ}). We may then introduce the sets of test functions:
\bea
\overline{\cA}u(t,\m) &:=& \Big\{\f \in C_b^{1,2}(\ol \bQ_t):  (\f-u)(t,\m) = \max_{\cN_\d(t, \m)}(\f-u) \ \mbox{for some $\d > 0$} \Big\} \nonumber, \\ 
\underline{\cA}u(t,\m) &:=& \Big\{\f \in C_b^{1,2}(\ol \bQ_t):  (\f-u)(t,\m) = \min_{\cN_\d(t, \m)}(\f-u) \ \mbox{for some $\d > 0$} \Big\}. \nonumber
\eea

\no This notion of viscosity solution enjoys the following useful properties: \\
$\bullet$ It is tailor-made for path-dependent PDEs. Although we chose to present first the argument for Markovian PDEs, we emphasize that, by using the same notion of viscosity solutions, both cases can be handled in a unified approach and with only notational modifications and minor changes in the proofs. \\
$\bullet$ It allows to easily construct semi-jets and strict extrema on the viscosity neighborhood, which are crucial to the proofs of our main convergence result. \\
$\bullet$ We can show (see Proposition \ref{lem-propagation}) that the points lying in the ``viscosity neighborhood" of the finite dimensional semi-solutions converge (up to a subsequence) to points in the viscosity neighborhood of semi-solutions of the mean field equations. \\

\no Although it could be possible to resort to other notions of viscosity solutions, we believe that this one is particularly adapted to our methodology. Requiring global tangency of test functions on the Wasserstein space (as in \cite{CGKPR}) would make difficult the use of semi-jets (see Proposition \ref{prop-semijets}), which greatly simplify the proof of our main result. The notion of \cite{BIRS} allows to use semi-jets but would require an extension to the case of path-dependent PDEs to be able to provide a unified approach for both Markovian and non-Markovian PDEs.

\begin{defn}
\label{defn-viscosity}
Let $u : \ol \bQ_0 \rightarrow \dbR$.\\
{\rm (i)} $u$ is a viscosity supersolution of \reff{pde} if, for all $(t, \m) \in \bQ_0$ and $\f \in \overline{\cA}u(t,m)$,
\beaa
\label{super}
 -\pa_t \f(t,\m) - F\big(t, \m, u(t,\m), \pa_x \d_m \f(t,\m,\cdot), \pa_{xx}^2 \d_m \f(t,\m,\cdot) \big) \ge 0. 
 \eeaa
  {\rm (ii)} $u$ is a viscosity subsolution of \reff{pde} if, for all $(t, \m) \in \bQ_0$ and $\f \in \underline{\cA}u(t,m)$,
\beaa
\label{sub}
 -\pa_t \f(t,\m) - F\big(t, \m, u(t,\m), \pa_x \d_m \f(t,\m,\cdot), \pa_{xx}^2 \d_m \f(t,\m,\cdot) \big) \le 0.  
 \eeaa
{\rm (iii)} $u$ is a viscosity solution of \eqref{pde} if it is a viscosity supersolution and subsolution. 

\end{defn}

\paragraph{Definition of viscosity solutions via semi-jets.} Fix $(t, \m) \in \bQ_0$ and $\d > 0$. For $(v, a, f) \in \dbR \times \dbR \times C_b^2(\dbR^d)$, introduce 
\bea\label{semijets}
\psi^{v,a, f}(s, \n) := v + a(s-t) + \langle \n - \m, f \rangle \q \mbox{for all $(s, \n) \in \cN_\d(t,\m)$.}
\eea
We then have the equivalence result: 
\begin{prop}\label{prop-semijets}
Let $u : \bQ_0 \longrightarrow \dbR$.\\
\no{\rm (i)} $u$ is a viscosity supersolution of \eqref{pde} if and only if it satisfies the viscosity supersolution property for all test functions in $\cup_{(t, \m) \in \bQ_0} \ol \cA u(t,\m)$ of the form \eqref{semijets}.\\
\no{\rm (ii)} $u$ is a viscosity subsolution of \eqref{pde} if and only if it satisfies the viscosity subsolution property for all test functions in  $\cup_{(t, \m) \in \bQ_0} \ul \cA u(t,\m)$  of the form \eqref{semijets}.
\end{prop}

\proof
We only provide the argument for (i). If $u$ is a viscosity supersolution of \eqref{pde}, then it satisfies the supersolution property for all $\f \in \ol \cA u(t,\m)$, in particular for those of the form \eqref{semijets}.

Assume now that the supersolution property is verified for all $\psi^{v, a, f} \in \ol \cA u(t,\m)$, $(v, a, f) \in \dbR \times \dbR \times C_b^2(\dbR^d)$. Let $\f \in \ol \cA u(t,\m)$ and $\e > 0$, and fix $(v, a, f) := \big(\f(t,\m), \pa_t \f(t, \m) - \e, \d_m \f(t, \m, \cdot)\big)$. We have, for $(s, \n) \in \cN_\d(t,\m)$:
\beaa
(\f - \psi^{v,a,f})(s, \n) &=& \f(s, \n) - \f(t, \m) - a(s-t) - \langle \n - \m, f \rangle \\
&=& \big(\f(s, \n) - \f(t, \n) - a(s-t)\big) + \big(\f(t, \n) - \f(t, \m) - \langle \n - \m, f \rangle \big) \\
&=& (s-t)\big(\pa_t \f(t, \n) + \eta(s-t) - \pa_t \f(t, \m) + \e \big)  \\ && + \int_0^1 \big\langle \n - \m, \d_m \f(t, \l \m + (1-\l) \n, \cdot) - \d_m \f(t, \m, \cdot) \big\rangle d\l
\eeaa
where $\eta(s-t) \underset{s \to t}{\longrightarrow} 0$. As $(s, \n) \in \cN_\d(t,\m)$, there exists $\dbP \in \cP(t,\m)$ s.t. $\n = \dbP_{X_s}$. Thus, introducing $h_{t,s}^\l := \d_m \f(t, \l \m + (1-\l) \n, \cdot) - \d_m \f(t, \m, \cdot)$, we have
$$\big\langle \n - \m, \d_m \f(t, \l \m + (1-\l) \n, \cdot) - \d_m \f(t, \m, \cdot) \big\rangle = \dbE^\dbP\Big[ h_{t,s}^\l(X_s) -  h_{t,s}^\l(X_t)  \Big]. $$ 
As $\f$ is smooth, we may apply Itô's formula to $h_{t,s}^\l$, and thus
\beaa
\dbE^\dbP\Big[  h_{t,s}^\l(X_s) -  h_{t,s}^\l(X_t) \Big] = \dbE^\dbP\Big[  \int_t^s \pa_x h_{t,s}^\l (X_r) \cdot dX_r + \frac 1 2  \pa_{xx}^2 h_{t,s}^\l (X_r) : d\langle X \rangle_r  \Big] \ge -(s-t) \frac \e 2,
\eeaa
for all $s \in [t, t+\d]$ and $\d$ sufficiently small, given the boundedness of the characteristics of $X$ under $\dbP$, the boundedness and continuity of the derivatives of $\f$ and the continuity of the flow $s \mapsto \dbP_{X_s}$. Finally, as we also have $\pa_t \f(t, \n) + \eta(s-t) - \pa_t \f(t, \m) \ge -\frac \e 2$ for $\d$ sufficiently small, we have 
$$ (\f - \psi^{v,a,f})(s, \n) \ge (s-t)\big( - \frac \e 2 + \e - \frac \e 2  \big) = 0 \q \mbox{for all $(s, \n) \in \cN_\d(t,\m),$} $$
which implies since $\psi^{v,a,f}(t, \m) = \f(t,\m)$ that $\psi^{v,a,f} \in \ol \cA u(t,\m)$. Then, the supersolution property writes
$$ - (\pa_t \f(t,\m)-\e) - F\big(t, \m, u(t,\m), \pa_x \d_m \f(t,\m, \cdot), \pa_{xx}^2 \d_m \f(t,\m, \cdot) \big) \ge 0, $$
and we obtain the desired result by letting $\e \longrightarrow 0$. 
\qed 

\paragraph{Comparison principle}

Under additional assumptions on $F$ (see \cite[Assumption 3.1]{WZ}), viscosity solutions satisfy the usual properties of consistency with the classical solution and stability. However, there is (at our knowledge) no uniqueness result for the general equation \eqref{pde}. Wu \& Zhang proved it in our setting for some specific cases (see \cite[Theorem 4.13]{WZ}). Over the past few years, many efforts have been main to obtain more general comparison results, see e.g.\ Burzoni, Ignazio, Reppen \& Soner \cite{BIRS}, Soner \& Yan \cite{soner2023viscosity}, Bertucci \cite{bertucci2023stochastic}, Bayraktar, Ekren \& Zhang \cite{bayraktar2023comparison}, Daudin \& Seeger \cite{daudin2023comparison}. In this paper, we shall particularly refer to the recent work of Cosso, Gozzi, Kharroubi, Pham \& Rosestolato \cite{CGKPR}, who established the comparison principle for Hamilton-Jacobi-Bellman equations on Wasserstein space in a quite general framework, assuming continuity of the semi-solutions. We refer to Remark \ref{rem:comparison} below for more detail about how their result relates to our notion of viscosity solution.
  In the statement of our main results, we shall use the comparison principle as a standing assumption:


\bas[Comparison principle]\label{comparison}
Let $u, v$ be respectively continuous viscosity subsolution and supersolution of \eqref{pde} such that $u(T, \cdot) \le v(T, \cdot)$. Then $u \le v$.
\eas

\begin{rem}\label{rem:comparison}{\rm
In the setting of \cite{CGKPR}, the tangency property of the test functions consists in global maxima/minima on $\ol \bQ_t$; thus, since a global extremum is \textit{a fortiori} an extremum on $\cN_\d(t,\m)$, any viscosity subsolution (resp. supersolution) in the sense of Definition \ref{defn-viscosity} is a viscosity subsolution (resp. supersolution) in the sense of \cite{CGKPR}, as long as we allow test functions to have derivatives with quadratic growth in $x$ instead of bounded ones (our requirement for boundedness has for purpose to provide a unified approach for both Markovian and path-dependent cases; however, our convergence result still holds under the quadratic growth requirement in the Markovian case, see Remark \ref{rem:quadratic}). Therefore, under the assumptions of \cite[Theorem 5.1]{CGKPR}, comparison holds for our notion of viscosity solutions.} \\
\end{rem}

\section{Finite-dimensional approximation}\label{sect-approximation}

\subsection{The approximating equation}
Let $N \ge 1$. We shall write in bold character the elements $\bx = (x_1, \dots, x_N) \in \dbR^{d \times N}, \bz = (z_1, \dots, z_N)$ and $\bgamma = {\rm Diag}(\g_1, \dots, \g_N) \in \dbS_{d \times N}^{\rm D}$. Introduce $F^N : [0,T] \times \dbR^{d \times N} \times \dbR \times \dbR^{d \times N} \times \dbS_{d \times N}^{\rm D} \longrightarrow \dbR$ such that:
\bea\label{cvg-assumption}
F^N\Big(t', \bx, y', \frac{\f(\bx)}{N}, \frac{\f'(\bx)}{N}\Big)\longrightarrow F(t, \m, y, \f, \f') \ \mbox{as} \ (N,t', \m^N(\bx), y') \longrightarrow (+\infty, t, \m, y)
\eea
for all $\f \in C_b^1(\dbR^d, \dbR)$, where we denote: 
$$\f(\bx) := (\f(x_1), \dots, \f(x_N)) \q \mbox{and} \q f'(\bx) = {\rm Diag}(\f'(x_1), \dots, \f'(x_N)).$$
\paragraph{A natural approximation.} Let us explain how to construct an approximating operator satisfying the consistency requirement \eqref{cvg-assumption}. Introduce, for $(t, \bx, y, \bz, \bgamma) \in [0,T] \times \dbR^{d \times N} \times \dbR \times \dbR^{d \times N} \times \dbS_{d \times N}^{\rm D}$,
\bea\label{Noperator}
F^N(t, \bx, y, \bz, \bgamma) := F(t, \m^N(\bx), y, N \bz \cdot \1_\bx, N\bgamma \cdot \1_\bx),
\eea
where $\bz \cdot \1_\bx(x) := \sum_{k=1}^N z_k\1_{x_k}(x)$ and $\bgamma \cdot \1_\bx(x) := \sum_{k=1}^N \g_k\1_{x_k}(x)$ for all $x \in \dbR^d$. Then:
\begin{prop}
Let Assumption \ref{assum-operator} hold. Then $F^N$ defined in \eqref{assum-operator} satisfies \eqref{cvg-assumption}.
\end{prop}
\proof
Let $\f \in C_b^1(\dbR^d, \dbR)$. We have:
\beaa
F_N\Big(t, \bx, y,  \frac{\f(\bx)}{N}, \frac{\f'(\bx)}{N}\Big) &=& F\big(t, \m^N(\bx), y, \sum_{k=1}^N \f(x_k) \1_{x_k}, \sum_{k=1}^N \f'(x_k) \1_{x_k}\big) \\
&=& F(t, \m^N(\bx), y, \f, \f')
\eeaa
by Assumption \ref{assum-operator} (ii). Then \eqref{cvg-assumption} is comes from the continuity assumption \ref{assum-operator} (i).
\qed


\no We now introduce the PDE on $[0,T] \times \dbR^{d \times N}$:
\bea\label{Npde}
-\pa_t u(t,\bx) - F^N\big(t, \bx, u(t,\bx), \pa_\bx u(t,\bx), \pa_{\bx \bx}^2 u(t,\bx)\big) = 0, \ u|_{t=T} = g^N,
\eea
with $g^N(\bx) := g(\m^N(\bx))$, $\pa_\bx u(t,\bx) := (\pa_{x_1} u, \dots, \pa_{x_N} u)(t,\bx) \in \dbR^{d \times N}$ and $\pa_{\bx \bx}^2 u(t,\bx) := {\rm Diag}(\pa_{x_1 x_1}^2 u, \dots, \pa_{x_N x_N}^2 u)(t,\bx) \in \dbS_{d \times N}$.  

\subsection{Viscosity solutions}\label{sec-Nvisc}

We define viscosity solutions for the equation \eqref{Npde}, as in the non-Markovian PDEs. We refer to Ren, Touzi \& Zhang \cite{RTZ} for a general overview of viscosity solutions for such equations. Let $\bX := (X^1, \dots, X^N)^\top$ be the canonical process on $\O^N$ and $\dbF^N = \{\cF_t^N\}_{t \in [0,T]}$ the corresponding filtration. For $t \in [0,T)$, define
$$ \L_t^N := [t,T) \times \dbR^{d \times N}, \q \mbox{and} \q \bar \L_t^N := [t,T] \times \dbR^{d \times N}. $$
\begin{defn}
For $(t, \bx) \in \L_0^N$, let $\cP_L^N(t, \bx)$ be the set of $\dbP \in \cP_2( \O^N, \cF_T^N)$ such that \\
$\bullet$ $\bX_t = \bx$, $\dbP$-a.s., \\
$\bullet$ there exist $(b^\dbP, \si^\dbP) :  [0,T] \times \O^N \longrightarrow \dbR^{d \times N} \times \dbS_{d \times N}^{\rm D}$, $\dbF^N$-measurable, bounded by $L$ coordinate-wisely, s.t.
\bea\label{cPLN}
d\bX_s =  b^\dbP_s ds +  \si^\dbP_s dW^\dbP_s,
\eea
where $W^\dbP$ is a $d \times N$-dimensional $\dbP$-Brownian motion. 
\end{defn}

\begin{lem}\label{lem-Ncompact}
The set $\cP_L^N(t, \bx)$ is weakly compact.
\end{lem}
\proof
Let $\tilde \cP_L^N(t, \bx)$ be defined as $\cP_L^N(t, \bx)$, without requiring that $\si^\dbP$ is blockwise diagonal. Clearly $\cP_L^N(t, \bx) \subset \tilde \cP_L^N(t, \bx)$, and we know from Zheng \cite[Theorem 3]{ZhengTightness} that $\tilde \cP_L^N(t, \bx)$ is weakly compact. Therefore, we only need to prove that $\cP_L^N(t, \bx)$ is closed under the weak convergence.

Let $(\dbP^n)_{n \ge 1}$ be sequence in $\cP_L^N(t, \bx)$ converging weakly to some $\dbP$, and denote $t \mapsto VT_t(Y)$ the total variation process associated with a process $Y$. Clearly, the family $\{ \dbP^n \circ \big(VT_t(\int_{0}^{.}  b_s^{\dbP^n}ds)\big)^{-1}\}_{n \ge 1}$ is tight for all $t \in [0,T]$ as the $b^{\dbP^n}$ are uniformly bounded. Therefore, we may apply Jacod \& Shiryaev \cite[Theorem 6.26]{jacod2013limit} to deduce that $\dbP^n_{\langle \bX \rangle}$ converges weakly to $\dbP_{\langle \bX \rangle}$. This implies in particular that $\si^\dbP$ still takes its values in $\dbS_{d \times N}^{\rm D}$, and therefore that $\cP_L^N(t, \bx)$ is closed under the weak convergence. \qed


\no Let $\cT_{t,T}^N$ denote the set of $[t,T]$-valued $\dbF^N$-stopping times, and $\cT_{t,T}^{N,+} :=  \{H \in \cT^N_{t,T} : H > t \}$. 
We define the sets of test functions:
\bea
\overline{\cA}^N  u(t,\bx) \!\!\!\!\! &:=& \!\!\!\!\! \Big\{\phi \in C_b^{1,2}(\ol \L_t^N): \exists H \in \cT_{t,T}^{N,+} \ \mbox{s.t.} \ (\phi-u)(t,\bx) = \max_{\th \in \cT^N_{t,T}} \ol \cE_{t,\bx}^N\Big[(\phi-u)(\th \wedge H, \bX_{\th \wedge H})\Big]\Big\} \nonumber, \\ 
\underline{\cA}^Nu(t,\bx) \!\!\!\!\! &:=& \!\!\!\!\! \Big\{\phi \in C_b^{1,2}(\ol \L_t^N): \exists H \in \cT_{t,T}^{N,+} \ \mbox{s.t.} \  (\phi-u)(t,\bx) =  \min_{\th \in \cT^N_{t,T}} \ul \cE_{t,\bx}^N\Big[(\phi-u)(\th \wedge H, \bX_{\th \wedge H })\Big]\Big\}, \nonumber
\eea
where $\ol \cE^{N}$ and $\ul \cE^{N}$ are the nonlinear expectations defined by
\bea\label{nonlinearexp}
\ol \cE^N_{t,\bx}[\cdot] := \sup_{\dbP \in \cP_L^N(t,\bx)} \dbE^\dbP[\cdot], \q \ul \cE^N_{t,\bx}[\cdot] := \inf_{\dbP \in \cP_L^N(t,\bx)} \dbE^\dbP[\cdot],
\eea
and $C_b^{1,2}(\ol \L_t^N)$ denotes the set of bounded functions of $C^{1,2}(\ol \L_t^N)$ with bounded derivatives. 

\begin{defn}
\label{defn-viscosityN}
Let $u : \bar \L_0^N \rightarrow \dbR$. \\
{\rm (i)} $u$ is a viscosity supersolution of \reff{Npde} if, for all $(t, \bx) \in \L_0^N$ and $\phi \in \overline{\cA}^Nu(t,\bx)$,
\bea
\label{super1}
-\pa_t \phi(t,\bx) - F^N\big(t, \bx, \phi(t,\bx), \pa_\bx \phi(t,\bx), \pa_{\bx \bx}^2 \phi(t,\bx)\big) \ge 0. 
 \eea
 {\rm (ii)} $u$ is a viscosity subsolution of \reff{Npde} if, for all $(t, \bx) \in \L_0^N$ and $\phi \in \underline{\cA}^Nu(t,\bx)$,
\bea
\label{sub1}
-\pa_t \phi(t,\bx) - F^N\big(t, \bx, \phi(t,\bx), \pa_\bx \phi(t,\bx), \pa_{\bx \bx}^2 \phi(t,\bx)\big) \le 0. 
 \eea
{\rm (iii)} $u$ is a viscosity solution of \eqref{Npde} if it is a viscosity supersolution and subsolution.
\end{defn}

\begin{rem}{\rm
The reader might find surprising our choice to resort to test functions for path-dependent PDEs in the context of Markovian equations. The reason is the following: this notion only involves points of the space writing as $\bX_{\th \wedge H}$ for some stopping time $\th$. Then, we prove in Proposition \ref{lem-propagation} that the sequence $\{\m^N(\bX_{\th \wedge H})\}_{N \ge 1}$ is tight, and that all its accumulation points lie in a (infinite dimensional) viscosity neighborhood $\cN_\d(t, \m)$. This property is crucial in the proof of our main results in Section \ref{sect-proof}.} \qed
\end{rem}

\subsection{Main results}

Let $\cS^N$ be the set of functions $h : \bar \L_0^N \longrightarrow \dbR$ such that $h(t,\bx) = {\rm h}^N(t, \m^N(\bx))$ for some ${\rm h}^N : [0,T] \times \cP^N(\dbR^d) \longrightarrow \dbR$. 
\begin{defn}\label{def-locunibounded}
We say $\{h^N\}_{N \ge 1} \in \prod_{N \ge 1} \cS^N$ is locally uniformly bounded if, for all $(t, \m) \in \ol \bQ_0$, there exist $\d,M > 0$ s.t., for all $(s, \bx^N) \in [0,T] \times \dbR^{d \times N}$ s.t. $\lvert s - t \rvert + \cW_2(\m^N(\bx^N), \m) \le \d$ and $N \ge 1$, we have $\lvert h^N(s, \bx^N) \rvert \le M$. 
\end{defn}

\no We now state the main result of the paper:

\begin{thm}\label{sol-cv}
Let $\{V^N \in \cS^N\}_{N \ge 1}$ be a sequence of 
continuous and locally uniformly bounded viscosity solutions of \eqref{Npde} s.t. $V^N|_{t=T} = g^N$, and introduce for all $(t,\m) \in \ol \bQ_0$
\bea\label{semilim}
\ul V(t,\m) := \underset{\tiny \begin{array}{c} N \to \infty, s \to t \\  \m^N(\bx^N) \overset{\cW_2}{\longrightarrow} \m \end{array}}{\lim \inf} {\rm V}^N(s, \m^N(\bx^N)), \q  \ol V(t,\m) := \underset{\tiny \begin{array}{c} N \to \infty, s \to t \\ \m^N(\bx^N) \overset{\cW_2}{\longrightarrow} \m \end{array}}{\lim \sup} {\rm V}^N(s, \m^N(\bx^N)). 
\eea
 If Assumption \ref{comparison} holds,  $\ul V$ and $\ol V$ are continuous and $\ul V|_{t=T} = \ol V|_{t=T} = g$,
then $V^N$ converges to the unique continuous viscosity solution $V$ of \eqref{pde}, i.e., the following limit exists:
$$  V(t,\m) = \underset{\tiny \begin{array}{c} N \to \infty, s \to t \\ \m^N(\bx^N) \overset{\cW_2}{\longrightarrow} \m \end{array}}{\lim} {\rm V}^N(s, \m^N(\bx^N)) \q \mbox{for all $(t,\m) \in \ol \bQ_0$} $$
and is the unique viscosity solution of \eqref{pde}. 
\end{thm}


\no The proof of this Theorem relies heavily on the following result, which corresponds to an adaptation of the Barles \& Souganidis \cite{BarSou} monotone scheme to our context:

\begin{thm}\label{visc-cv}
{\rm (i)} Let $\{v^N \in \cS^N\}_{N \ge 1}$ be a sequence of continuous and locally uniformly bounded 
viscosity supersolutions of \eqref{Npde}.
The relaxed semi-limit defined by  
$$\ul v(t,\m) := \underset{\tiny \begin{array}{c} N \to \infty, s \to t \\ \m^N(\bx^N) \overset{\cW_2}{\longrightarrow} \m \end{array}}{\lim \inf} {\rm v}^N(s, \m^N(\bx^N)) \q \mbox{for all $(t,\m) \in \ol \bQ_0$}$$ 
is finite and is a LSC 
viscosity supersolution of $\eqref{pde}$. \\
{\rm (ii)} Let $\{u^N \in \cS^N\}_{N \ge 1}$ be a sequence of continuous and locally uniformly bounded 
viscosity subsolutions of \eqref{Npde}. 
The relaxed semi-limit defined by  
$$\ol u(t,\m) := \underset{\tiny \begin{array}{c} N \to \infty, s \to t \\ \m^N(\bx^N) \overset{\cW_2}{\longrightarrow} \m \end{array}}{\lim \sup}  {\rm u}^N(s, \m^N(\bx^N)) \q \mbox{for all $(t,\m) \in \ol \bQ_0$}$$ 
is finite and is a USC 
viscosity subsolution of $\eqref{pde}$. 
\end{thm}
The proof of this result is relegated to Section \ref{sect-proof}.


\begin{rem}[Comparison with the Barles-Souganidis monotone scheme]{\rm
(i) Our finite-dimensional approximation shares strong similarities with the numerical scheme of \cite{BarSou} for second order PDEs. The condition \eqref{cvg-assumption} can be seen as the {\it consistency} condition, and the existence of locally uniformly bounded solutions to \eqref{Npde} as the {\it stability} condition. However, the {\it monotonicity} condition seems less obvious at first sight. Note that our ``scheme" is defined as viscosity solutions to a PDE, rather than as a classical solution to an approximating equation such as (2.1) in \cite{BarSou}. We then believe that the monotonicity condition lies in the very fact that $V^N$, for $N \ge 1$, is a viscosity solution to \eqref{Npde}; thus, for test functions tangent to $V^N$ from above, we have the inequality \eqref{sub1}, and the converse inequality \eqref{super1} for test functions tangent from below. 

\no (ii) The main motivation of the scheme of \cite{BarSou} is to derive numerical approximations for PDEs. It is then natural to wonder how the present result could be used to achieve this objective, i.e.\ finding numerical approximations of the solution of the PDE on Wasserstein space \eqref{pde}. Since our approximating function $V^N$ is defined as a viscosity solution to a (finite-dimensional) second order PDE, one natural idea is the following: $V^N$ can be approximated by the monotone scheme of \cite{BarSou}. The numerical approximation is then a function $V_\rho^N$, $\rho > 0$, with $V_\rho^N \to V^N$ as $\rho \to 0$. The question then boils down to finding an efficient numerical scheme to approximate $V^N$ for large $N$. This could legitimately be attempted via deep learning methods, see e.g.\ \cite{sirignano2018dgm} or \cite{carmona2021convergence}. 
\qed
}
\end{rem}

\section{Application to stochastic control}\label{sect-examples}

\subsection{Mean field control}

 Let $k \ge 1$, $A \subset \dbR^k$ and $(b,\si): [0,T] \times \dbR^d \times \cP_2(\dbR^d) \times A \rightarrow \dbR^d, \dbS_d$, continuous in $(t,a) \in [0,T] \times A$, Lipschitz-continuous in $(x,\m) \in \dbR^d \times \cP_2(\dbR^d)$ uniformly in $(t,a)$ and uniformly bounded by $L$. For $(t,\m) \in \bQ_0$ and $\a : [0,T] \times \O \longrightarrow A$, let $\dbP^{t, \m, \a}$ be s.t. $X$ is a controlled McKean-Vlasov diffusion with drift and diffusion coefficients $b$ and $\si$, i.e.
\bea\label{controlledMKV}
X_s = \xi + \int_t^s b(r, X_r, \dbP^{t, \m, \a}_{X_r}, \a_r)dr + \int_t^s \si(r, X_r, \dbP^{t, \m, \a}_{X_r}, \a_r)dW_r^\a, \ \mbox{$\dbP^{t, \m, \a}$-a.s,} 
\eea
where $W^\a$ is a standard $d$-dimensional $\dbP^{t,\m,\a}$-Brownian motion and $\dbP_\xi^{t, \m, \a} = \m$. Let $\cA_t$ be the set of $\dbF$-progressively measurable processes $\a : [t,T] \times \O \longrightarrow A$ such that \eqref{controlledMKV} has a unique weak solution.
 We consider the mean field control problem
\bea\label{MFcontrol}
V(t, \m) := \sup_{\a \in \cA_t} \dbE^{\dbP^{t, \m, \a}}\Big[ \int_t^T f(r, X_r, \dbP^{t, \m, \a}_{X_r}, \a_r)dr + g(\dbP^{t, \m, \a}_{X_T}) \Big], \nonumber 
\eea 
with $f : [0,T] \times \dbR^d \times \cP_2(\dbR^d) \times A \longrightarrow \dbR$ and $g : \cP_2(\dbR^d) \longrightarrow \dbR$. We know from Wu \& Zhang \cite[Theorem 5.8]{WZ} that, if $V$ is continuous, then it is a viscosity solution of the following Hamilton-Jacobi-Bellman (HJB) equation on Wasserstein space:
\bea\label{MFHJB}
-\pa_t u(t, \m)-F_{\rm HJB}(t, \m, \pa_x \d_m u(t, \m, \cdot), \pa_{xx}^2 \d_m u(t, \m, \cdot)) = 0, \ u(T, \cdot) = g,
\eea
where
\bea\label{FHJB}
F_{\rm HJB}(t, \m, Z, \G) := \Big\langle \m, \ \sup_{a \in A} \ \Big\{ b(t, \cdot, \m, a) \cdot Z(\cdot) + \frac{1}{2}\si^2(t, \cdot, \m, a):\G(\cdot) + f(t, \cdot, \m, a) \Big\} \Big\rangle.
\eea
\begin{prop}\label{prop-assumF}
Assume $A$ is compact. Then $F_{\rm HJB}$ satisfies Assumption \ref{assum-operator}.
\end{prop}
\proof 
Let $(t^n, \m^n, Z^n, \G^n)$ be a sequence converging to some $(t, \m, Z, \G) \in \bQ_0 \times C^0(\dbR^d, \dbR^d) \times C^0(\dbR^d, \dbS^d)$ in the sense of Assumption \ref{assum-operator} (i). Observe that:
\beaa
&&\Big\lvert \int_{\dbR^d} \sup_{a \in A} \ b(t^n, x, \m^n, a)Z^n(x)\m^n(dx) - \int_{\dbR^d} \sup_{a \in A} \ b(t, x, \m, a)Z(x)\m(dx) \Big\rvert \\ &&\le \Big\lvert \int_{\dbR^d} \sup_{a \in A} \ b(t^n, x, \m^n, a)Z^n(x)(\m^n - \m)(dx) \Big\rvert + \Big\lvert \int_{\dbR^d} \sup_{a \in A} \ b(t^n, x, \m^n, a)(Z^n - Z)(x)\m(dx) \Big\rvert \\ && \q + \Big\lvert \int_{\dbR^d} \sup_{a \in A} \ \big(b(t^n, x, \m^n, a) - b(t, x, \m, a)\big)Z(x)\m(dx) \Big\rvert \\
&& \le C\int_{\dbR^d} (1 + x^2)(\m^n - \m)(dx) + L \int_{\dbR^d} \lvert Z^n(x) - Z(x) \rvert \m(dx)  \\ && \q + \Big\lvert \int_{\dbR^d} \sup_{a \in A} \ \big(b(t^n, x, \m^n, a) - b(t, x, \m, a)\big)Z(x)\m(dx) \Big\rvert,
\eeaa
for some constant $C$ independent from $n$. The first term of the right-hand side converges to $0$ because $\cW_2(\m^n, \m) \to 0$; the second term converges to $0$ by the dominated convergence theorem, as $Z^n$ converges pointwise to $Z$ and the functions $\{Z^n\}_{n \ge 0}$ and $Z$ have quadratic growth uniformly in $n$; finally, the third term converges to $0$ because it is continuous in $(t^n, \m^n)$, by continuity of $b$ and compactness of $A$. Using similar estimates to handle terms in $\si$ and $f$, we deduce that $F_{\rm HJB}$ satisfies Assumption \ref{assum-operator} (i). As to (ii), it is clearly satisfied as $F_{\rm HJB}$ is an integral w.r.t. $\m$.
\qed

%

\begin{rem}{\rm
We observe that, in the case of $F_{\rm HJB}$, $Z$ and $\G$ may belong to $\dbL^1(\m)$. However, we chose to restrict them to sets of bounded functions when we introduced the operator $F$ in order to have a more general framework and avoid possible integrability issues.} \qed
\end{rem}

\paragraph{Finite-dimensional approximation}

For $(t,\bx) \in [0,T] \times \dbR^{d \times N}$ and a given control $\balpha = (\a^1, \dots, \a^N) : [0,T] \times \O^N \longrightarrow A^N$, let $\dbP^{t, \bx, \balpha}$ be such that, for all $i \in [N] := \{1, \dots, N\}$, 
\bea\label{NcontrolledMKV}
X_s^i = x_i + \int_t^s b(r, X_r^i, \m^N(\bX_r), \a_r^i)dr + \int_t^s \si(r, X_r^i, \m^N(\bX_r), \a_r^i)dW_r^{i,\balpha}, \ \mbox{$\dbP^{t, \bx, \balpha}$-a.s,} 
\eea
where $W^{\balpha} := (W^{1, \balpha}, \dots, W^{N, \balpha})^\top$ is a standard $d \times N$-dimensional $\dbP^{t, \bx, \balpha}$-Brownian motion. Let $\cA_t^N$ be the set of $\dbF^N$-progressively measurable processes $\balpha : [t,T] \times \O^N \longrightarrow A^N$ s.t. \eqref{NcontrolledMKV} has a unique weak solution. We define the control problem
\bea\label{MFcontrolN}
V^N(t, \bx) := \sup_{\balpha \in \cA_t^N} \sum_{i=1}^N \dbE^{\dbP^{t, \bx, \balpha}}\Big[ \int_t^T f^{i,N}(r, \bX_r, \a_r^i)dr + g^N(\bX_T) \Big], \nonumber 
\eea 
with $f^{i,N}(t, \bx, a) := f(r, x_i, \m^N(\bx), a)$.
We know from standard stochastic control theory that, if $V^N$ is continuous, then it is a viscosity solution of
\bea\label{NHJB}
- \pa_t u(t,\bx) - \!\! \sup_{\ba \in A^N}\Big\{ {\bf b}(t, \bx, \ba) \cdot \pa_{\bx} u(t, \bx) + \frac 1 2 \bsigma^2(t, \bx, \ba) : \pa_{\bx \bx}^2 u(t,\bx) + {\bf f}(t, \bx, \ba) \cdot \be \Big\} = 0,
\eea
 where
  ${\bf b}(t, \bx, \ba) := \Big(b\big(t, x_i, \m^N(\bx), a_i\big)\Big)_{1 \le i \le N}^\top$, $\bsigma(t, \bx, \ba) := {\rm Diag}\Big(\si \big(t, x_i, \m^N(\bx), a_i\big)\Big)_{1 \le i \le N}$, ${\bf f}(t,\bx, \ba) := \Big(f^{i,N}\big(t, \bx, a_i\big)\Big)_{1 \le i \le N}^\top$ and $\be := (1, \dots, 1)^\top \in \dbR^N$.
  

 \begin{prop}\label{prop-cvcontrol}
 Assume that:\\
 $\bullet$ $f$ and $g$ are bounded and continuous on $\cP_2(\dbR^d)$, and extend continuously on $\cP_1(\dbR^d)$;\\ 
 $\bullet$ $b$, $\si$ and $f$ are $\b$-Hölderian in $t$, uniformly in $(x, \m, a) \in \dbR^d \times \cP_2(\dbR^d) \times A$, for some $\b \in (0,1]$;\\
 $\bullet$ $\si$ does not depend on $\m$, and satisfies $\si(\cdot, a) \in C^{1,2}([0,T] \times \dbR^d)$ for all $a \in A$, with all its derivatives uniformly bounded.\\
 Then $V^N$ converges to $V$, i.e.,
 $$ V(t, \m) =   \underset{\tiny \begin{array}{c} N \to \infty, s \to t \\ \cW_2(\m^N(\bx^N), \m) \to 0 \end{array}}{\lim}  V^N(s,\bx^N) \q \mbox{for all $(t,\m) \in \bQ_0$}. $$
 \end{prop}
 \proof
 We first show that \eqref{NHJB} is the finite-dimensional approximation of \eqref{MFHJB}, i.e.,
$$ F_{\rm HJB}^N(t, \bx, \bz, \bgamma) = \sup_{\ba \in A^N}\Big\{ {\bf b}(t, \bx, \ba) \cdot \bz + \frac 1 2 \bsigma^2(t, \bx, \ba) : \bgamma + {\bf f}(t, \bx, \ba) \cdot \be \Big\}, $$
where $F_{\rm HJB}^N$ is the finite-dimensional approximation of $F_{\rm HJB}$ defined by \eqref{Noperator}. We compute
\beaa
F_{\rm HJB}^N(t, \bx, \bz, \bgamma) &=& F_{\rm HJB}(t, \m^N(\bx), N\bz \cdot \1_\bx, N\bgamma \cdot \1_\bx) \\
&& \hspace{-3cm} = \Big\langle \m^N(\bx), \ N \sup_{a \in A} \ \Big\{ b(t, \cdot, \m^N(\bx), a) \cdot \bz \cdot \1_\bx+ \frac{1}{2}\si^2(t, \cdot, \m^N(\bx), a): \bgamma \cdot \1_\bx + f(t, \cdot, \m^N(\bx), a) \Big\} \Big\rangle \\
&& \hspace{-3cm} = \sum_{i=1}^N  \sup_{a \in A} \Big\{ b(t, x_i, \m^N(\bx), a) \cdot z_i+ \frac{1}{2}\si^2(t, x_i, \m^N(\bx), a) : \g_i  + f(t, x_i, \m^N(\bx), a) \Big\} \\
&& \hspace{-3cm} = \sup_{{\bf a} \in A^N} \sum_{i=1}^N \Big\{ b(t, x_i, \m^N(\bx), a_i) \cdot z_i+ \frac{1}{2}\si^2(t, x_i, \m^N(\bx), a_i) : \g_i + f(t, x_i, \m^N(\bx), a_i) \Big\} \\
&& \hspace{-3cm} = \sup_{\ba \in A^N}\Big\{ {\bf b}(t, \bx, \ba) \cdot \bz + \frac 1 2 \bsigma^2(t, \bx, \ba) : \bgamma + {\bf f}(t, \bx, \ba) \cdot \be \Big\}.
\eeaa
Moreover, as $f$ and $g$ are bounded,
 the $V^N$ are uniformly bounded, and $\ul V$ and $\ol V$ (defined similarly to \eqref{semilim}) are bounded.

We now prove that $\ul V$ and $\ol V$ are continuous. Let $(\O^0, \cF^0, \dbP^0)$ be a probability space such that we can construct for all $(t, \bx, \balpha) \times [0,T] \times \dbR^{d \times N} \times \cA_t^N$ a diffusion process $\bX^{t, \bx, \balpha}$ such that $\dbP^0 \circ (\bX^{t, \bx, \balpha})^{-1} = \dbP^{t, \bx, \balpha} \circ \bX^{-1}$.
Fix $\bx$, $\bx'$ and $R > 0$ s.t.\ $\lVert \bx \rVert_2, \lVert \bx' \rVert_2 < R$, and $R' > R$ to be determined later. Observe that $D_{R'} := \{ m \in \cP_2(\dbR^d) : \lVert m \rVert_2 \le R' \}$ is bounded in $\cP_2(\dbR^d)$, and therefore is $\cW_1$-compact. Thus, there exists a continuity modulus $\rho_{R'}$ for $g$ on this set, and then:
\beaa 
&&\dbE^{\dbP^0}\Big[ \big\lvert g\big( \m^N(\bX_T^{t, \bx, \balpha}) \big) - g(\m^N(\bX_T^{t, \bx', \balpha}))  \big\rvert \Big] \\ &\le& \dbE^{\dbP^0}\Big[ \rho_{R'}\big(\cW_1\big( \m^N(\bX_T^{t, \bx, \balpha}), \m^N(\bX_T^{t, \bx', \balpha}) \big)\big) \\&&+ \big\lvert g\big( \m^N(\bX_T^{t, \bx, \balpha}) \big) - g(\m^N(\bX_T^{t, \bx', \balpha}))  \big\rvert \big(\1_{D_{R'}^c}(\m^N(\bX_T^{t, \bx, \balpha})) + \1_{D_{R'}^c}(\m^N(\bX_T^{t, \bx', \balpha}))\big) \Big]  \\
&\le& \dbE^{\dbP^{0}}\Big[ \rho_{R'}\big(\cW_1\big( \m^N(\bX_T^{t,\bx,\balpha}), \m^N(\bX_T^{t,\bx',\balpha}) \big)\big) + C \dbP^0\big( \lVert \bX_T^{t,\bx,\balpha} \rVert_2 \ge R') + C\dbP^0\big(\lVert \bX_T^{t,\bx,\balpha} \rVert_2 \ge R') \Big]  \\
&\le&  \dbE^{\dbP^{0}}\Big[ \rho_{R'}\big(\cW_1\big( \m^N(\bX_T^{t,\bx,\balpha}), \m^N(\bX_T^{t,\bx',\balpha}) \big)\big) \Big]+ \frac{C'}{R'}(1 + R^2) 
\eeaa
with $C' >0$ is independent from $N$ and $\balpha$, due to the uniform boundedness of the drift and diffusion coefficients, and to Markov's inequality. Note also that
\beaa
&& \dbE^{\dbP^{0}}\Big[ \rho_{R'}\big(\cW_1\big( \m^N(\bX_T^{t,\bx,\balpha}), \m^N(\bX_T^{t,\bx',\balpha}) \big)\big) \Big] \\ && \qq \le  \rho_{R'}(\eta) + \frac{1}{\sqrt{\eta}}\sqrt{ \dbE^{\dbP^{0}}\Big[ \rho_{R'}^2\big(\cW_1\big( \m^N(\bX_T^{t, \bx, \balpha}), \m^N(\bX_T^{t, \bx', \balpha})  \big)\big)\Big] \dbE^{\dbP^{0}}\Big[\cW_1^2\big( \m^N(\bX_T^{t, \bx, \balpha}), \m^N(\bX_T^{t, \bx', \balpha}) \big)  \Big] }
\eeaa
 for all $\eta > 0$. Fix $\e > 0$. Choosing $R' = R_\e := \frac{2C'(1+R^2)}{\e}$, we can find $\d > 0$ and $\eta > 0$ (possibly depending on $\e$ but not on $N$) such that
 $$ \dbE^{\dbP^{0}}\Big[ \big\lvert g\big( \m^N(\bX_T^{t,\bx,\balpha}) \big) - g(\m^N(\bX_T^{t,\bx',\balpha}))  \big\rvert \Big]  \Big\rvert \le \e \q \mbox{whenever $\lVert \bx - \bx' \rVert_2 \le \d$}. $$
 We may prove a similar estimate with $f$, and finally, by arbitrariness of $\balpha \in \cA^N$:
\bea\label{estimate-loc-unif}
 \Big\lvert V^N(t, \bx) - V^N(t, \bx') \Big\rvert \le \e \q \mbox{whenever $\lVert \bx - \bx' \rVert_2 \le \d$}. \nonumber
 \eea
 Using similar estimates, we may also prove that  
 \bea\label{estimate-loc-unif1}
 \Big\lvert V^N(t, \bx) - V^N(t', \bx) \Big\rvert \le \e \q \mbox{whenever $|t - t'| \le \d$}, \nonumber
 \eea
 which implies that $\ul V$ and $\ol V$ are continuous and satisfy $\ul V|_{t=T} = \ol V|_{t=T} = g$. Note also that, by symmetry of the problem, $V^N \in \cS^N$. Since we are under the assumptions of \cite[Theorem 5.1]{CGKPR}, by Remark \ref{rem:comparison}, Assumption \ref{comparison} holds. We may then apply Theorem \ref{sol-cv} to derive the convergence result.  \qed

 \begin{rem}{\rm
 We emphasize that the assumptions on $b$ and $\si$ and the Hölder regularity assumption on $f$ in Proposition \ref{prop-cvcontrol} are solely used to apply the comparison principle of \cite{CGKPR}. Therefore, any availability of a less restrictive comparison principle in the literature would automatically relax the assumptions we need to ensure this convergence result. \qed
 }
 \end{rem}

 \subsection{Zero-sum stochastic differential games}
 
 We present in a more informal way a second example. We consider the following control problem, arising in zero-sum games:
 $$ V_+(t, \m) := \inf_{\a^1 \in \cA_t^1} \sup_{\a^2 \in \cA_t^2} \dbE^{\dbP^{t, \m, \a^1, \a^2}}\Big[ \int_t^T f\big(s, X_s, \dbP_{X_s}^{t, \m, \a^1, \a^2}, \a_s^1,\a_s^2\big)ds + g\big(\dbP_{X_T}^{t, \m, \a^1, \a^2}\big)\Big],  $$
 where the measures $\{ \dbP^{t, \m, \a^1, \a^2} : \a_1 \in \cA_t^1, \a_2 \in \cA_t^2\}$ are such that $X$ has the dynamics \eqref{controlledMKV}, substituting $(\a^1, \a^2)$ to $\a$.
 By Cosso \& Pham \cite{CosPha}, $V_+$ is a viscosity solution of \eqref{pde}, with operator
\bea\label{Fplus}
F_+(t, \m, Z, \G) := \Big\langle \m, \inf_{a_2 \in A_2} \sup_{a_1 \in A_1} \!\!\! \Big\{ b(t, \cdot, \m, a_1, a_2) \cdot Z + \frac{1}{2}\si^2(t, \cdot, \m, a_1, a_2):\G + f(t, \cdot, \m, a_1, a_2) \Big\} \Big\rangle. \nonumber
\eea
Although \cite{CosPha} uses another notion of viscosity solutions, we may consider in the context of our discussion that $V_+$ is also a viscosity solution in the sense of Definition \ref{defn-viscosity}. Assuming $A_1$ and $A_2$ are compact, the corresponding finite-dimensional approximation is then given by \eqref{Npde}, with operator
 $$ F_{+}^N(t, \bx, \bz, \bgamma) := \inf_{\ba_1 \in A_1^N} \sup_{\ba_1 \in A_2^N}\Big\{ {\bf b}(t, \bx, \ba_1, \ba_2) \cdot \bz + \frac 1 2 \bsigma^2(t, \bx, \ba_1, \ba_2) : \bgamma + {\bf f}(t, \bx, \ba_1, \ba_2) \cdot \be \Big\}. $$
 As in the case of mean field control, we may show that this corresponds to the PDE satisfied by the control problem
$$ V_+^N(t, \bx) := \inf_{\balpha^1 \in (\cA_t^1)^N} \sup_{\balpha^2 \in (\cA_t^2)^N} \sum_{i=1}^N \dbE^{\dbP^{t, \bx, \balpha^1, \balpha^2}}\Big[ \int_t^T f^{i,N}\big(r, \bX_r, \balpha_r^1,\balpha_r^2\big)dr + g\big(\bX_T\big)\Big],  $$
 where the measures $\{ \dbP^{t, \bx, \balpha^1, \balpha^2} : \balpha_1 \in \cA_t^1, \balpha_2 \in \cA_t^2\}$ are such that $\bX$ has the dynamics \eqref{NcontrolledMKV}, substituting $(\balpha^1, \balpha^2)$ to $\a$.
 
  \begin{prop}
 If Assumption \ref{comparison} holds for \eqref{pde} with $F = F_+$, then $V_+^N$ converges to $V_+$, i.e.,
 $$ V_+(t, \m) =   \underset{\tiny \begin{array}{c} N \to \infty, s \to t \\ \cW_2(\m^N(\bx^N), \m) \to 0 \end{array}}{\lim}  V_+^N(s,\bx^N) \q \mbox{for all $(t,\m) \in \bQ_0$}. $$
 \end{prop}

\subsection{The uncontrolled case}

The purpose of this paragraph is to recover the classical propagation of chaos result for diffusion processes, whose first instance was given by Snitzman \cite{Snitzman} for some special models. $b$ and $\si$ do no longer depend on the variable $a$. We consider the equation
\bea\label{linear}
-\pa_t u(t,\m) \! - \! \Big\langle \m,  b(t, \cdot, \m) \! \cdot \! \pa_x \d_m u(t,\m,\cdot) + \frac{1}{2}\si^2(t, \cdot, \m) \! : \! \pa_{xx}^2 \d_m u(t, \m,\cdot)  \Big\rangle = 0, \ u(T, \cdot) = g,
\eea
where $g \in C_b^0(\cP_2(\dbR^d), \dbR)$, the set of continuous and bounded functions from $\cP_2(\dbR^d)$ to $\dbR$. For $N \ge 1$, we easily see that the corresponding finite-dimensional approximation writes
\bea\label{Nlinear}
-\pa_t u^N(t,\bx) -  {\bf b}(t, \bx) \cdot \pa_\bx u^N(t,\bx) - \frac{1}{2}{\bsigma}^2(t, \bx):\pa_{\bx \bx}^2 u^N(t, \bx) = 0, \ u^N(T, \cdot) = g^N.
\eea
For $(t,\m, \bx^N) \in \ol \bQ_0 \times \dbR^{d \times N}$, let $(\bar \dbP^{t,\m}, \bar \dbP^{t,\bx^N}) \in \cP_L(t,\m) \times \cP_L^N(t,\bx^N)$ be such that $X$ and $\bX$ are the uncontrolled versions of \eqref{controlledMKV} and \eqref{NcontrolledMKV}, respectively $\bar \dbP^{t,\m}$-a.s. and $\bar \dbP^{t,\bx^N}$-a.s. As $g \in C_b^0(\cP_2(\dbR^d), \dbR)$, we know that, under some smoothness assumptions on $b$ and $\si$ (see Talbi, Touzi \& Zhang \cite[Lemma 3.6]{TTZ2}), we have
$$ u(t,\m) = g\big(\bar \dbP^{t,\m}_{X_T}\big), \ u^N(t, \bx^N) = \dbE^{\bar \dbP^{t,\bx^N}}\big[ g^N(\bX_T) \big] = \dbE^{\bar \dbP^{t,\bx^N}}\big[ g\big(\m^N(\bX_T^N)\big) \big], $$
for all $(t,\m) \in \ol \bQ_0$ and $\bx^N \in \dbR^{d \times N}$. Thus, applying Proposition \ref{prop-cvcontrol}, we have $u^N(0, \bx^N) \longrightarrow u(0, \m)$ as $N \to \infty$ and $\cW_2(\m^N(\bx^N), \m) \to 0$, hence 
$$ \dbE^{\bar \dbP^{0,\bx^N}}\big[ g\big(\m^N(\bX_T^N)\big) \big] \underset{N \to \infty}{\longrightarrow} g\big(\bar \dbP^{0,\m}_{X_T}\big), $$
which exactly means that $\bar \dbP^{0,\bx^N} \circ (\m^N(\bX_T^N))^{-1}$ converges weakly to $\bar \dbP^{0,\m}_{X_T}$, as it is true for all $g \in C_b^0(\cP_2(\dbR^d), \dbR)$. This corresponds to the propagation of chaos result proved by Oelschlager \cite{Oelschlager}. 

\section{Extension to path-dependent PDEs}\label{sect-PPDE}

\subsection{Pathwise derivatives}

For $t \in [0,T)$, we adapt our previous notations to the path-dependent case: 
$$ \bQ_t := [t,T) \times \cP_2(\O) \q \mbox{and} \q \ol \bQ_t := [t,T] \times \cP_2(\O). $$
For $(t,\m) \in \ol \bQ_0$, we denote by $\m_{[0,t]}$ the law of the stopped process $X_{ \cdot \wedge t}$ under $\m$. We shall use the notion of pathwise derivative of Ekren, Keller, Touzi \& Zhang \cite{EKTZ}, which is tailor-made for continuous semimartingales. In particular, it allows to introduce a notion of derivative that is intrinsic to the space of continuous paths, whereas the notion of Dupire \cite{Dupire} which requires to include càdlàg paths.

\begin{defn}\label{def-pathderiv}
{\rm (i)} Given a metric space $E$, we denote by $C^0([0,T] \times \O, E)$ the set of $\dbF$-progressively measurable and continuous functions from $[0,T] \times \O$ to $E$, where $\O$ is equipped with the norm $\lvert \o \rvert := \sup_{t \in [0,T]} \lvert \o_t \rvert$ for all $\o \in \O$.\\
{\rm (ii)} We denote by $u \in C^{1,2}([0,T] \times \O)$ the set of functions $u : [0,T] \times \O \to \dbR$ such that there exist $\pa_t u \in C^0([0,T] \times \O, \dbR)$, $\pa_\o u \in C^0([0,T] \times \O, \dbR^d)$ and $\pa_{\o \o}^2 u \in C^0([0,T] \times \O, \dbS_d)$ such that, for all $\dbP \in \bigcup_{L > 0} \cP_L$, $u$ satisfies
$$ du(t,X) = \pa_t u(t,X)dt + \pa_\o u(t,X) \cdot dX_t + \frac 1 2 \pa_{\o \o}^2 u(t,X) : d\langle X \rangle_t, \q \mbox{$\dbP$-a.s.} $$ 
\end{defn}

 \begin{defn}\label{C12Sp}
 Fix $t \in [0,T]$.
 
\no{\rm (i)} We denote by $C^0(\ol \bQ_t)$ the set of functions $u : \ol \bQ_t \to \dbR$ continuous for the pseudo-metric:
\bea\label{pseudometric}
 \widetilde \cW_2\big((s, \m), (r, \n)\big) := \Big( \lvert s - r \rvert^2 + \cW_2^2(\m_{[0,s]}, \n_{[0,r]}) \Big)^{\frac 1 2} \q \mbox{for all $(s, \m),(r,\n) \in \bQ_t$.}
\eea
{\rm (ii)} We denote by $C^{1,2}_b(\ol \bQ_t)$ the set of bounded functions $u: \ol \bQ_t \to \dbR$ such that $\pa_t u$, $\d_m u$, $\pa_\o \d_m u, \pa_{\o\o}^2\d_m u$ exist, are bounded in all their variables and continuous in $(t,\m)$ in the sense of {\rm (i)}, where the functional linear derivative takes the form $\d_m u: [t, T]\times \cP_2(\O)\times \O \to \dbR$ satisfying, for any $s\in [t, T]$ and $\m, \m' \in \cP_2(\O)$,
\beaa
u(s,\m')-u(s, \m) = \displaystyle\int_0^1 \int_{\O} \d_m u(s, \l \m' + (1-\l)\m, \o)(\m'-\m)(d\o)d\l.
\eeaa
\end{defn}
Note that, if $u \in C^0(\ol \bQ_t)$, then $u(s, \m) = u(s, \m_{[0,s]})$ for all $(s, \m) \in \ol \bQ_t$.

\subsection{Path-dependent equation on Wasserstein space}

Let $F : [0,T] \times \cP_2(\O) \times \dbR \times \dbL_2^0(\O, \dbR^d) \times \dbL_2^0(\O, \dbS_d) \longrightarrow \dbR$, where $\dbL_2^0(\O,\dbR^d)$ (resp. $\dbL_2^0(\O, \dbS_d)$) denotes the set of Borel measurable functions from $\O$ to $\dbR^d$ (resp. $\dbS_d$) with quadratic growth,  and $g : \cP_2(\O) \longrightarrow \dbR$.  We consider the following equation:
\bea\label{ppde}
-\pa_t u(t,\m) - F\big(t, \m, u(t,\m), \pa_\o \d_m u(t,\m,\cdot), \pa_{\o\o}^2 \d_m u(t,\m,\cdot) \big) = 0, \ u|_{t=T} = g, \ \mbox{$(t,\m) \in \bQ_0$.} 
\eea
We define semijets similarly to \eqref{semijets} and straightforwardly adapt Proposition \ref{prop-semijets} and Assumption \ref{assum-operator} to the path-dependent setting.

\subsection{Viscosity solutions}

We redefine, for all $(t,\m) \in \bQ_0$, $\cP_L(t,\m) := \{ \dbP \in \cP_L : \dbP_{X_{t \wedge \cdot}} = \m_{[0,t]} \}$, as well as the neighborhood
$$ \cN_\d (t,\m) := [t,t+\d] \times \cP_L(t,\m), $$
which is compact under $\widetilde \cW_2$ (see again Wu \& Zhang \cite[Lemma 4.1]{WZ}). We then introduce the sets of test functions:
\bea
\overline{\cA}u(t,\m) &:=& \Big\{\f \in C_b^{1,2}(\ol \bQ_t):   (\f-u)(t,\m) = \max_{\cN_\d(t, \m)}(\f-u) \ \mbox{for some $\d > 0$}\Big\} \nonumber, \\ 
\underline{\cA}u(t,\m) &:=& \Big\{\f \in C_b^{1,2}(\ol \bQ_t):   (\f-u)(t,\m) = \min_{\cN_\d(t, \m)}(\f-u) \ \mbox{for some $\d > 0$}\Big\}. \nonumber
\eea

\begin{defn}
\label{defn-viscosityp}
Let $u : \ol \bQ_0 \rightarrow \dbR$.\\
{\rm (i)} $u$ is a viscosity supersolution of \reff{ppde} if, for all $(t, \m) \in \bQ_0$ and $\f \in \overline{\cA}u(t,\m)$,
\beaa
\label{super}
 -\pa_t \f(t,\m) - F\big(t, \m, u(t,\m), \pa_\o \d_m \f(t,\m,\cdot), \pa_{\o\o}^2 \d_m \f(t,\m,\cdot) \big) \ge 0. 
 \eeaa
 {\rm (ii)} $u$ is a viscosity subsolution of \reff{ppde} if, for all $(t, \m) \in \bQ_0$ and $\f \in \underline{\cA}u(t,\m)$,
\beaa
\label{sub}
 -\pa_t \f(t,\m) - F\big(t, \m, u(t,\m), \pa_\o \d_m \f(t,\m,\cdot), \pa_{\o\o}^2 \d_m \f(t,\m,\cdot) \big) \le 0.  
 \eeaa
 {\rm (iii)} $u$ is a viscosity solution of \eqref{ppde} if it is a viscosity supersolution and subsolution.

\end{defn}
\subsection{Finite-dimensional approximation}

Let $N \ge 1$. We shall write in bold character the elements $\bo = (\o^1, \dots, \o^N) \in \O^N$.  Introduce $F^N : [0,T] \times \O^N \times \dbR \times \dbR^{d \times N} \times \dbS_{d \times N}^{\rm D} \longrightarrow \dbR$ such that:
\bea\label{cvg-assumption-p}
F^N\Big(t', \bo, y', \frac{\f(\bo)}{N}, \frac{\f'(\bo)}{N}\Big)\longrightarrow F(t, \m, y, \f, \f') \ \mbox{as} \ (N,t', \m^N(\bo), y') \longrightarrow (+\infty, t, \m, y)
\eea
for all $\f \in C_b^1(\dbR^d, \dbR)$, where we denote: 
$$\f(\bo) := (\f(\o^1), \dots, \f(\o^N)) \q \mbox{and} \q f'(\bx) = \rm{Diag}(\f'(\o^1), \dots, \f'(\o^N)).$$ 
As in the Markovian case, we may guarantee the existence of an approximation. Introduce, for $(t, \bo, y, \bz, \bgamma) \in [0,T] \times \O^N \times \dbR \times \dbR^{d \times N} \times \dbS_{d \times N}^{\rm D}$,
\bea\label{Noperatorp}
F^N(t, \bo, y, \bz, \bgamma) := F(t, \m^N(\bo), y, N\bz \cdot \1_\bo, N\bgamma \cdot \1_\bo),
\eea
where $\bz \cdot \1_\bo(\o) := \sum_{k=1}^N z_k\1_{\o^k}(\o)$ and $\bgamma \cdot \1_\bo(\o) := \sum_{k=1}^N \g_k\1_{\o^k}(\o)$ for all $\o \in \O$.  Then, if Assumption \ref{assum-operator} holds, then $F^N$ satisfies \ref{cvg-assumption-p}.

\no We now introduce the path-dependent PDE on $[0,T] \times \O^N$:
\bea\label{Nppde}
-\pa_t u(t,\bo) - F^N\big(t, \bo, u(t,\bo), \pa_\bo u(t,\bo), \pa_{\bo\bo}^2 u(t,\bo)\big) = 0, \ u|_{t=T} = g^N,
\eea
with $g^N(\bo) := g(\m^N(\bo))$, $\pa_\bo u(t,\bo) := (\pa_{\o^1} u, \dots, \pa_{\o^N} u)(t,\bo) \in \dbR^{d \times N}$ and $\pa_{\bo \bo}^2 u(t,\bo) := {\rm Diag}(\pa_{\o^1 \o^1}^2 u, \dots, \pa_{\o^N \o^N}^2 u)(t,\bo) \in \dbS_{d \times N}$. 

We now define viscosity solutions for \eqref{Nppde}. We adapt the notations of Section \ref{sec-Nvisc} to the path-dependent case.  For $t \in [0,T)$, define
$$ \L_t^N := [t,T) \times \O^N, \q \mbox{and} \q \bar \L_t^N := [t,T] \times \O^N. $$
For $(t, \bo) \in \L_0^N$, we define $\cP_L^N(t, \bo)$ similarly to $\cP_L^N(t,\bx)$, with the condition $\bX_{t \wedge \cdot} = \bo_{t \wedge \cdot}$, $\dbP$-a.s,, for all $\dbP \in \cP_L^N(t,\bo)$.
Similarly to Lemma \ref{lem-Ncompact}, we have the following result:
\begin{lem}\label{lem-Ncompactppde}
The set $\cP_L^N(t, \bo)$ is weakly compact.
\end{lem}
We define the sets of test functions:
\bea
\overline{\cA}^N  u(t,\bo) \!\!\!\!\! &:=& \!\!\!\!\! \Big\{\phi \in C_b^{1,2}(\ol \L_0^N): \exists H \in \cT_{t,T}^{N,+} \ \mbox{s.t.} \ (\phi-u)(t,\bo) = \max_{\th \in \cT^N_{t,T}} \ol \cE_{t,\bo}^N\Big[(\phi-u)(\th \wedge H, \bX_{\cdot \wedge \th \wedge H})\Big]\Big\} \nonumber, \\ 
\underline{\cA}^Nu(t,\bo) \!\!\!\!\! &:=& \!\!\!\!\! \Big\{\phi \in C_b^{1,2}(\ol \L_0^N): \exists H \in \cT_{t,T}^{N,+} \ \mbox{s.t.} \  (\phi-u)(t,\bo) =  \min_{\th \in \cT^N_{t,T}} \ul \cE_{t,\bo}^N\Big[(\phi-u)(\th \wedge H, \bX_{\cdot \wedge \th \wedge H})\Big]\Big\}, \nonumber
\eea
where $\ol \cE^{N}$ and $\ul \cE^{N}$ are the nonlinear expectations defined by
\bea\label{nonlinearexpp}
\ol \cE^N_{t,\bo}[\cdot] := \sup_{\dbP \in \cP_L^N(t,\bo)} \dbE^\dbP[\cdot], \q \ul \cE^N_{t,\bo}[\cdot] := \inf_{\dbP \in \cP_L^N(t,\bo)} \dbE^\dbP[\cdot],
\eea
and $C_b^{1,2}(\ol \L_t^N)$ denotes the bounded elements of $C^{1,2}(\ol \L_t^N)$ (defined similarly to $C^{1,2}([0,T] \times \O)$) with bounded derivatives.

\begin{defn}
\label{defn-viscosityN}
Let $u : \bar \L_0^N \rightarrow \dbR$. \\
{\rm (i)} $u$ is a viscosity supersolution of \reff{pde} if, for all $(t, \bo) \in \L_0^N$ and $\phi \in \overline{\cA}^Nu(t,\bo)$,
\bea
-\pa_t \phi(t,\bo) - F^N\big(t, \bo, \phi(t,\bo), \pa_\bo \phi(t,\bo), \pa_{\bo \bo}^2 \phi(t,\bo)\big) \ge 0. \nonumber
 \eea
 {\rm (ii)} $u$ is a viscosity subsolution of \reff{pde} if, for all $(t, \bo) \in \L_0^N$ and $\phi \in \underline{\cA}^Nu(t,\bo)$, 
\bea
-\pa_t \phi(t,\bo) - F^N\big(t, \bo, \phi(t,\bo), \pa_\bo \phi(t,\bo), \pa_{\bo \bo}^2 \phi(t,\bo)\big) \le 0. \nonumber
 \eea
{\rm (iii)} $u$ is a viscosity solution of \eqref{pde} if it is a viscosity supersolution and subsolution.

\end{defn}

In this paragraph, $\cS^N$ denotes the set of functions $h : \bar \L_0^N \longrightarrow \dbR$ s.t. $h(t,\bo) = {\rm h}^N(t, \m^N(\bo))$ for some ${\rm h}^N : [0,T] \times \cP^N(\O) \longrightarrow \dbR$. 

\begin{thm}\label{sol-cvp}
Let $\{V^N \in \cS^N\}_{N \ge 1}$ be a sequence of uniformly
continuous for \eqref{pseudometric} and locally bounded, uniformly in $N$, viscosity solutions of \eqref{Nppde} s.t. $V^N(T, \cdot) = g^N$, and introduce
\bea\label{semilim2}
\ul V(t,\m) := \underset{\tiny \begin{array}{c} N \to \infty, s \to t \\  \m^N(\bo^N) \overset{\cW_2}{\longrightarrow} \m \end{array}}{\lim \inf} {\rm V}^N(s, \m^N(\bo^N)), \q  \ol V(t,\m) := \underset{\tiny \begin{array}{c} N \to \infty, s \to t \\ \m^N(\bo^N) \overset{\cW_2}{\longrightarrow} \m \end{array}}{\lim \sup} {\rm V}^N(s,  \m^N(\bo^N)). 
\eea
 If Assumption \ref{comparison} holds and $\ul V|_{t=T} = \ol V|_{t=T} = g$,
then $u^N$ converges to the unique continuous viscosity solution $V$ of \eqref{ppde}, i.e.,
$$  V(t,\m) = \underset{\tiny \begin{array}{c} N \to \infty, s \to t \\ \m^N(\bo^N) \overset{\cW_2}{\longrightarrow} \m \end{array}}{\lim} {\rm V}^N(s, \m^N(\bo^N)) \q \mbox{for all $(t,\m) \in \bQ_0$}. $$ 
\end{thm}

\begin{thm}\label{visc-cvp}
{\rm (i)} Let $\{v^N \in \cS^N\}_{N \ge 1}$ be a sequence of uniformly continuous for \eqref{pseudometric} and locally bounded, uniformly in $N$, 
viscosity supersolutions of \eqref{Nppde}. 
Then, the relaxed semi-limit defined by  
$$\ul v(t,\m) := \underset{\tiny \begin{array}{c} N \to \infty, s \to t \\ \m^N(\bo^N) \overset{\cW_2}{\longrightarrow} \m \end{array}}{\lim \inf} {\rm v}^N(s,  \m^N(\bo^N)) \q \mbox{for all $(t,\m) \in \bQ_0$}$$ 
is finite and is a LSC 
viscosity supersolution of $\eqref{pde}$. \\
{\rm (ii)} Let $\{ u^N \in \cS^N \}_{N \ge 1}$ be a sequence of uniformly continuous for \eqref{pseudometric} and locally bounded, uniformly in $N$, 
viscosity subsolutions of \eqref{Nppde}. 
Then, the relaxed semi-limit defined by  
$$\ol u(t,\m) := \underset{\tiny \begin{array}{c} N \to \infty, s \to t \\  \m^N(\bo^N) \overset{\cW_2}{\longrightarrow} \m \end{array}}{\lim \sup} {\rm u}^N(s,  \m^N(\bo^N)) \q \mbox{for all $(t,\m) \in \bQ_0$}$$ 
is finite and is a USC 
viscosity subsolution of $\eqref{pde}$. 
\end{thm}
These results are proved in Section \ref{sect-proof}. 

\section{Proof of the main results}\label{sect-proof}

\subsection{The Markovian setting}

\no{\bf Proof of Theorem \ref{sol-cv}} (given Theorem \ref{visc-cv})
By Theorem \ref{visc-cv}, $\ul V$ and $\ol V$ are respectively continuous viscosity supersolution and subsolution of \eqref{pde}. As $\ul V|_{t=T} = \ol V|_{t=T} = g$, the comparison principle implies $\ul V \ge \ol V$. By \eqref{semilim}, we also have the converse inequality, and thus $\ul V = \ol V$, and these functions are viscosity solutions of \eqref{pde}. Given the comparison principle, \eqref{pde} has a unique continuous viscosity solution, and thus $\ul V = \ol V = V$ and the two semi-limits \eqref{semilim} are equal to the limit.
\qed

\no{\bf Proof of Theorem \ref{visc-cv}}
We only prove the convergence of the viscosity supersolutions, as the case of the subsolutions is handled similarly. It is clear that $\ul v$ is finite as $\{v^N\}_{N \ge 1}$ is locally bounded, uniformly in $N$. Fix $(t,\m) \in \bQ_0$ and $\f \in \ol \cA \ul v(t,\m)$ with corresponding $\d_0 \in (0,T-t)$. By Proposition \ref{prop-semijets}, we may assume w.l.o.g. that $\f$ is a semijet of the form \eqref{semijets}, with characteristics $(v, a, f) \in \dbR \times \dbR \times C_b^2(\dbR^d)$. We also introduce, for all $N \ge 1$, the functions $\phi^N(s, \bx) := \f(s, \m^N(\bx))$ for all $(s,\bx) \in \ol \L_0^N$. Finally, let $(t^N, \bx^N)$ be a sequence such that $t^N \longrightarrow t$, $\m^N(\bx^N) \overset{\cW_2}{\longrightarrow} \m$, and ${\rm v}^N(t^N, \m^N(\bx^N)) \longrightarrow \ul v(t,\m)$ as $N \to \infty$. 

The main idea is the following: we shall approximate the $\{t^N, \bx^N\}_{N \ge 1}$ with ``good" points in which we may apply the viscosity supersolution property of $v^N$, and then deduce the one of $\ul v$ by passing to the limit as $N \to \infty$ . Such good points are given by the following lemma:

\begin{lem}\label{lem-delta}
There exists a family $\{t_\d^N, \bx^{\d,N}\}$, $\d >0$ and $N \ge 1$, such that $\phi^N \in \ol \cA v^N(t_\d^N, \bx^{\d,N})$ and $(t_\d^N, \m^N(\bx^{\d,N})) \underset{\d \to 0}{\longrightarrow} (t^N, \bx^N)$ for all $N \ge 1$.
\end{lem}
Let $(t_\d^N, \bx^{\d, N})$ be as in the above lemma. The supersolution property of $v^N$ provides:
\bea\label{supercv1}
-\pa_t \phi^N(q_\d^N) - F^N\big(q_\d^N, v^N(q_\d^N),  \pa_\bx\phi^N(q_\d^N), \pa_{\bx \bx}^2 \phi^N(q_\d^N) \big)  \ge 0.
\eea
where $q_\d^N := (t_\d^N, \bx^{\d,N})$. By the equalities \eqref{derivatives1} below, we have:
\begin{align*}
&F^N\big(q_\d^N, v^N(q_\d^N),  \pa_\bx\phi^N(q_\d^N), \pa_{\bx \bx}^2 \phi^N(q_\d^N) \big) = F^N\Big(q_\d^N, v^N(q_\d^N), \frac{f'(\bx^{\d,N})}{N}, \frac{f''(\bx^{\d,N})}{N}\Big)  
\end{align*}
Note that  
$$(t_\d^N, \m^N(\bx^{\d,N}), {\rm v}^N(t_\d^N, \m^N(\bx^{\d,N}))) \underset{\d \to 0}{\longrightarrow} (t^N, \m^N(\bx^N),{\rm v}^N(t^N, \m^N(\bx^N)))$$ by Lemma \ref{lem-delta} and continuity of ${\rm v}^N$, and thus 
$$(t_\d^N, \m^N (\bx^{\d,N}),  {\rm v}^N(t_\d^N, \m^N(\bx^{\d,N}))) \underset{(\d, N) \to (0, \infty)}{\longrightarrow} (t,\m, \ul v(t,\m)).$$
 We then deduce from the consistency property \eqref{cvg-assumption} that
\begin{align*}
F^N\big(q_\d^N, v^N(q_\d^N),  \pa_\bx\phi^N(q_\d^N), \pa_{\bx \bx}^2 \phi^N(q_\d^N) \big) \underset{\tiny \begin{array}{c} \d \to 0 \\ N \to \infty \end{array}}{\longrightarrow} F\big( t, \m, \ul v(t,\m),  f', f'' \big).
\end{align*}
Finally, as $\pa_t \phi^N(q_\d^N) = a = \pa_t \f(t,\m)$,  sending $(\d, N) \longrightarrow (0,\infty)$ in \eqref{supercv1} provides the viscosity supersolution property of $\ul v$.  \qed

\no{\bf Proof of Lemma \ref{lem-delta}}
 Replacing $\f$ with $\tilde \f(s,\cdot) := \f(s, \cdot) - (s-t)^2$, we may also assume w.l.o.g. that $(t,\m)$ is a strict maximum of $(\f - \ul v)$ on $\cN_{\d_0}(t,\m)$. 

\no{\bf Step 1:} Fix $\d \in (0, \frac{\d_0}{3})$, and introduce the stopping time
\begin{align*}
H_\d^N := \inf \big\{ s \ge t^N : \cW_2\big(\m^N(\bX_s), \m^N(\bx^N) \big) = 2\d  \big\} \wedge (t^N + 2\d).
\end{align*}
By Lemma \ref{lem:convex} and pathwise continuity of $\bX$, we may choose $\d$ sufficiently small so that the domain of $H_\d^N$ is a convex set (namely, the Euclidian ball centered in $\bx^N$ with radius $2\d$). Then, by Ekren, Touzi \& Zhang \cite[Theorem 3.5]{ETZ3} and weak compactness of $\cP_L^N(t, \bx^N)$,
there exists $(\th_\d^N, \dbP^{N,*}) \in \cT^N_{t^N, T} \times \cP_L^N(t, \bx^N)$ s.t.
\bea\label{optPtheta}
\dbE^{\dbP^{N,*}}\Big[(\phi^N - v^N)(\th_\d^N  \wedge H_\d^N, \bX_{\th_\d^N  \wedge H_\d^N}) \Big] = \!\! \sup_{\th \in \cT^N_{t^N, T}} \!\! \ol \cE^N_{t^N,\bx^N}\Big[(\phi^N - v^N)(\th  \wedge H_\d^N, \bX_{\th  \wedge H_\d^N}) \Big],
\eea
 where $\ol \cE^N_{t^N,\bx^N}$ is defined by \eqref{nonlinearexp}. Indeed, by continuity of $\phi^N-v^N$, we easily see that the Markov process $s \mapsto (\phi^N-v^N)(s, \bX_s)$ is bounded and uniformly continuous on $\{(t, \bo) : t \le H_\d^N(\bo) \}$, and that $\cP_L^N(t, \bx^N)$ satisfies \cite[Assumption 3.4]{ETZ3}.  Also, note that, since $\{v^N\}_{N \ge 1}$ is locally bounded, uniformly in $N$, and $\m^N(\bx^N) \overset{\cW_2}{\longrightarrow} \m$, we may assume w.l.o.g. that, after passing to an appropriate subsequence and for $\d$ small enough, 
 $\{(\phi^N - v^N)(\th \wedge H_\d^N, \bX_{\th \wedge H_\d^N}) \}_{N \ge 1}$ is uniformly bounded for all $\th \in \cT^N_{t,T}$.

\no{\bf Step 2:} We now justify that
\bea\label{theta}
\underset{N \to \infty}{\lim \sup} \ \dbP^{N,*}\big( \th_\d^{N} <H_\d^N \big) > 0.
\eea
Indeed, assume to the contrary that $\underset{N \to \infty}{\lim \sup} \ \dbP^{N,*}\big( \th_\d^{N} <H_\d^N \big) = \underset{N \to \infty}{\lim} \ \dbP^{N,*}\big( \th_\d^{N} <H_\d^N \big) = 0$. 
We have:
\bea\label{estsup}
(\f - \ul v)(t,\m) &=& \underset{N \to \infty}{\lim} (\f - {\rm v}^N)(t^N, \m^N(\bx^N))  \nonumber \\ 
&\le& \underset{N \to \infty}{\lim \inf} \ \dbE^{\dbP^{N,*}}\Big[(\f - {\rm v}^N)(\th_\d^N   \wedge H_\d^N , \m^N(\bX_{\th_\d^N   \wedge H_\d^N}))\Big] \nonumber  \\
&\le& \underset{N \to \infty}{\lim \sup} \ \dbE^{\dbP^{N,*}}\Big[(\f - {\rm v}^N)(\th_\d^N  \wedge H_\d^N , \m^N(\bX_{\th_\d^N   \wedge H_\d^N }))\Big] \nonumber  \\
&=& \underset{N \to \infty}{\lim \sup} \ \dbE^{\dbP^{N,*}}\Big[\Big\{ (\f - {\rm v}^N)(H_\d^N, \m^N(\bX_{H_\d^N}))(1-\1_{\th_\d^N < H_\d^N}) \nonumber \\
&&\hspace{2.2cm}+ (\f - {\rm v}^N)(\th_\d^N, \m^N(\bX_{\th_\d^N}))\1_{\th_\d^N < H_\d^N} \Big\} \Big] \nonumber \\
&=&  \underset{N \to \infty}{\lim \sup} \ \dbE^{\dbP^{N,*}}\Big[ (\f - {\rm v}^N)( H_\d^N, \m^N(\bX_{H_\d^N}))\Big], \nonumber
\eea
where we used the fact that $\{(\f - {\rm v}^N)(\th_\d^N   \wedge H_\d^N, \m^N(\bX_{\th_\d^N   \wedge H_\d^N}))\}_{N \ge 1}$ is uniformly bounded 
and $\underset{N \to \infty}{\lim \sup} \ \dbP^{N,*}\big(\th_\d^{N} <H_\d^N \big) = \underset{N \to \infty}{\lim} \ \dbP^{N,*}\big(\th_\d^{N} <H_\d^N \big) = 0$. 

Since  $(t^N, \m^N(\bx^N)) \underset{N \to \infty}{\longrightarrow} (t,\m)$, by Proposition \ref{lem-propagation} and compactness of $[t,T]$, there exists a subsequence $\n^N := \dbP^{N,*} \circ \big( H_\d^N, \m^N(\bX) \big)^{-1}$ that converges weakly to some $\n \in \cP_2( [t,T] \times  \cP_2(\O))$ supported on $ [t, t+\d] \times \cP_L(t,\m)$.
Thus, denoting by $(\t, m)$ the canonical mapping on $[0,T] \times \cP_2(\O)$, we have by upper semicontinuity of $\f - v^N$ and continuity of $(\t, m) \mapsto m_{X_\t}$,
\bea\label{bar-omega}
\underset{N \to \infty}{\lim \sup} \ \dbE^{\dbP^{N,*}}\Big[ (\f - {\rm v}^N)( H_\d^N, \m^N(\bX_{H_\d^N}))\Big]  &=& \underset{N \to \infty}{\lim \sup} \ \dbE^{\n^N}\Big[ (\f - {\rm v}^N)(\t, m)\Big] \nonumber \\
 &&\hspace{-5cm} \le \dbE^{\n}\Big[ (\f - \ul v)(\t, m_{X_\t})\Big] 
 \le (\f - \ul v)(\t(\bar \o), m_{X_{\t(\bar \o)}}(\bar \o)), 
\eea
for some $\bar \o \in [0,T] \times \cP_2(\O)$. 
We also observe that $\t > t$, $\n$-a.s. Indeed, given that, by definition of $H_\d^N$, we have
\beaa
  \cW_2\big(m_{X_\t}, \m^N(\bx^N) \big) \vee (\t-t^N) = 2\d \ \mbox{or} \  \cW_2\big(m_{X_\t}, \m^N(\bx^N) \big) + (\t-t^N) \ge 2\d, \q \mbox{$\n^N$-a.s.,}
\eeaa
for all $N \ge 1$, and therefore, for $N$ sufficiently large,
\beaa
\cW_2\big(m_{X_\t}, \m \big) \vee (\t-t) \le 3\d < \d_0 \ \mbox{or} \ \cW_2\big(m_{X_\t}, \m \big) +(\t - t) \ge \d, \q \mbox{$\n^N$-a.s.}
\eeaa
These inequalities define a fixed closed support for $\n^N$, which is inherited by $\n$ by weak convergence and continuity of $\cW_2$ and $\o \mapsto m_{X_{\t(\o)}}(\o)$. Thus, $\bar \o$ in \eqref{bar-omega} may be chosen s.t.
$$ \cW_2\big(m_{X_{\t(\bar \o)}}(\bar \o), \m \big) \vee (\t(\bar \o)-t) \le 3\d < \d_0 \ \mbox{or} \ \cW_2\big(m_{X_{\t(\bar \o)}}(\bar \o), \m \big) +(\t(\bar \o) - t) \ge \d.$$
The first inequality shows that $(\t(\bar \o), m_{X_{\t(\bar \o)}}(\bar \o)) \in \cN_{\d_0}(t,\m)$, and the second one that $(\t(\bar \o), m_{X_{\t(\bar \o)}}(\bar \o)) \neq (t,m)$. Therefore, \eqref{bar-omega} contradicts the fact that $(t, \m)$ is a strict maximum on $\cN_{\d}(t,\m)$. Thus, \eqref{theta} holds true, and we may find a subsequence $\{\bo^{\d,N}\}_{N \ge 1}$ s.t. $t_\d^N := \th_\d^N(\bo^{\d,N}) <H_\d^N(\bo^{\d,N})$ for all $N \ge 1$. 


\no {\bf Step 3:} We now prove that $\phi^N$ is a test function for $v^N$ in some well chosen point. Introduce $\bx^{\d,N} := \bX_{\th_\d^N(\bo^{\d,N})}(\bo^{\d,N})$ and $\bY^N$, the nonlinear Snell envelop of $s \mapsto (\phi^N - v^N)(s, \bX_s)$, i.e.,
$$ \bY_s^N(\bo) := \sup_{\th \in \cT^N_{t,T}} \ol \cE_{s, \bo}^N\Big[ (\phi^N - v^N)(\th \wedge H_\d^N, \bX_{\th \wedge H_\d^N})\Big], $$
 which satisfies $\bY_s \ge \ol \cE_{s, \bY_s}^N[\bY_{\th \wedge H_\d^N}]$ for all $\th \in \cT^N_{s, t}$. Then we have, for all $\th \in \cT^N_{t_\d^N, T}$,
\beaa
(\phi^N - v^N)(t_\d^N, \bx^{\d,N}) &=& \bY_{t_\d^N}^N(\bo^{\d,N}) \\ &&\hspace{-4cm}\ge \ol \cE_{t_\d^N, \bo^{\d,N}}^N\big[ \bY_{\th  \wedge H_\d^N}^N \big] \ge \ol \cE_{t_\d^N, \bo^{\d,N}}^N\big[(\phi^N - v^N)(\th  \wedge H_\d^N, \bX_{\th  \wedge H_\d^N}) \big],
\eeaa
and therefore, as we are in the Markovian case,
$$(\phi^N - v^N)(t_\d^N, \bx^{\d,N}) = \max_{\th \in \cT^N_{t_\d^N,T}} \ol \cE^N_{t_\d^N, \bx^{\d,N}}\Big[(\phi^N - v^N)(\th  \wedge H_\d^N, \bX_{\th  \wedge H_\d^N}) \Big].  $$
Observe that $\phi^N \in C_b^{1,2}(\ol \L_{t_\d^N}^N)$. Indeed, since $\f = \psi^{v, a, f}$ is a semijet, we have 
\bea\label{derivatives1}
\pa_t \phi^N(s,\bx) &=& \pa_t \f(s, \m^N(\bx)) = a, \nonumber \\
\pa_{x_i} \phi^N(s,\bx) &=& \frac 1 N f'(x_i), \\ \nonumber 
\pa_{x_i x_i}^2 \phi^N(s, \bx)  &=& \frac 1 N  f''(x_i), 
\eea
 for all $i \in [N]$. As $H_\d^N > t_\d^N$ on $\{ \bX_{\cdot \wedge t_\d^N} = \bo_{\cdot \wedge t_\d^N}^{\d, N} \}$, we have $\phi^N \in \ol \cA^N v^N(t_\d^N, \bx^{\d,N})$. Finally, the fact that $(t_\d^N, \m(\bx^{\d, N})) \to (t^N, \bx^N)$ simply comes from the definition of $H_\d^N$. 
\qed

\begin{rem}\label{rem:quadratic}{\rm
If we allow our test functions on Wasserstein space to have derivatives with quadratic growth in $x$ (similarly to \cite{CGKPR}), we may not use our semijets, and therefore the computations in \eqref{derivatives} must be done differently. In the case of state-dependent functions, this can be done by using general formulas for smooth functions on the space of measures, see e.g.\ Carmona \& Delarue \cite[Vol. 1, Propositions 5.35 \& 5.91]{CarDel}. \qed}
\end{rem}

\subsection{The path-dependent setting}

\no{\bf Proof of Theorem \ref{sol-cvp}} (Given Theorem \ref{visc-cvp}). Identical to the proof of Theorem \ref{sol-cv}, using Theorem \ref{visc-cvp} instead of Theorem \ref{visc-cv}. \qed

\no{\bf Proof of Theorem \ref{visc-cvp}} This follows the same arguments as the proof of Theorem \ref{visc-cv}. We use the notations of Section \ref{sect-PPDE}. Let $(t^N, \bo^N)$ be a sequence s.t.\ $(t^N, \m^N(\bo^N)) \longrightarrow (t,\m)$ and ${\rm v}^N(t^N, \m^N(\bo^N)) \longrightarrow \ul v(t,\m)$ as $N \to \infty$. We also introduce, for all $N \ge 1$, the functions $\phi^N(s, \bo) := \f(s, \m^N(\bo))$ for all $(s,\bo) \in \L_0^N$. Fix $\d \in (0, \d_0)$, and define the stopping time
\begin{align*}
H_\d^N := \inf \big\{ s \ge t^N : \cW_2\big(\m^N(\bX_{\cdot \wedge s}), \m^N(\bo_{\cdot \wedge t^N}^N) \big) \ge \d  \big\} \wedge (t^N + \d).
\end{align*}
 Observe that, for all $s \ge t^N, \bo \in \O^N$ s.t. $s \le H_\d^N(\bo)$, we have
\beaa
 | v^N(s, \bo) | &\le& v^N(t^N, \bo^N) | + \rho_N\Big( |s - t^N| + \cW_2\big(\m^N(\bo_{\cdot \wedge s}), \m^N(\bo_{\cdot \wedge t^N}^N)\big)\Big) \\
 &\le& v^N(t^N, \bo^N) + \rho_N(2\d),
 \eeaa
 where $\rho_N$ is continuity modulus of $v^N$. Furthermore, $\phi^N$ is Lipschitz-continuous as it has bounded derivatives. Thus $\phi^N - v^N$ is bounded and uniformly continuous on $\{(s, \bo) : s \le H_\d^N(\bo)\}$, and by \cite{ETZ3} again,
there exists $(\th_\d^N, \dbP^{N,*}) \in \cT^N_{t^N, T} \times \cP_L^N(t, \bo^N)$ s.t.
\bea\label{optPthetap}
\dbE^{\dbP^{N,*}}\Big[(\phi^N - v^N)(\th_\d^N  \wedge H_\d^N, \bX_{\cdot \wedge \th_\d^N  \wedge H_\d^N}) \Big] = \!\! \sup_{\th \in \cT^N_{t^N, T}} \!\! \ol \cE^N_{t^N,\bo^N}\Big[(\phi^N - v^N)(\th  \wedge H_\d^N, \bX_{\cdot \wedge \th  \wedge H_\d^N }) \Big], \nonumber
\eea
 where $\ol \cE^N_{t^N,\bo^N}$ is defined by \eqref{nonlinearexpp}. 
Similarly to the Markovian setting, we may find a subsequence $\{\bo^{\d,N}\}_{N \ge 1}$ s.t. $t_\d^N := \th_\d^N(\bo^{\d,N}) <H_\d^N(\bo^{\d,N})$ for all $N \ge 1$, and satisfying 
$$(\phi^N - v^N)(t_\d^N, \bo^{\d,N}) = \max_{\th \in \cT^N_{t_\d^N, T}} \ol \cE^N_{t_\d^N, \bo^{\d,N}}\Big[(\phi^N - v^N)(\th  \wedge H_\d^N, \bX_{\cdot \wedge \th  \wedge H_\d^N}) \Big].  $$
Observe that $\phi^N \in C_b^{1,2}(\ol \L_{t_\d^N}^N)$. Indeed, since $\f = \psi^{v, a, f}$ is a semijet, we have 
\bea\label{derivatives}
\pa_t \phi^N(s,\bo) &=& \pa_t \f(s, \m^N(\bo)) = a, \nonumber \\
\pa_{\o_i} \phi^N(s,\bx) &=& \frac 1 N \pa_\o f(\o_i), \\ 
\pa_{\o_i \o_i}^2 \phi^N(s, \bo)  &=& \frac 1 N \pa_{\o\o}^2 f(\o_i), \nonumber 
\eea
 for all $i \in [N]$. As $H_\d^N > t_\d^N$ on $\{ \bX_{t_\d^N \wedge \cdot} = \bo_{t_\d^N \wedge \cdot}^{\d, N} \}$, we have $\phi^N \in \ol \cA^N v^N(t_\d^N, \bo^{\d,N})$ and the supersolution property provides
\bea\label{supercv}
-\pa_t \phi^N(q_\d^N) - F^N\big(q_\d^N, v^N(q_\d^N),  \pa_\bo\phi^N(q_\d^N), \pa_{\bo \bo}^2 \phi^N(q_\d^N) \big)  \ge 0.
\eea
where $q_\d^N := (t_\d^N, \bo^{\d,N})$. We conclude similarly to (i).
 \qed

\section{A precompactness result}\label{sect-propagation}

In this section, we state and prove our propagation of chaos-like result for continuous semimartingale with bounded characteristics, which plays a crucial role in the contradiction argument used to prove Lemma \ref{lem-delta}.

Our objective is to prove that the empirical measure associated with a $N$-dimensional continuous semimartingales with characteristics bounded by some constant $L$ converges in law (up to a subsequence) to an element supported on $\cP_L$, i.e.\ a measure on $\O$ under which the canonical process is also almost surely a continuous semimartingale with characteristics bounded by the same constant $L$. 

\begin{prop}\label{lem-propagation}
Let $\{\bo^N\}_{N \ge 1} \in \prod_{N \ge 1} \O^N$ and $\m \in \cP_2(\O)$ s.t. $\m^N(\bo_{\cdot \wedge t}^N) \overset{\cW_2}{\longrightarrow} \m_{[0,t]}$, and $\{\dbP^N \in \cP_L^N(t, \bo^N)\}_{N \ge 1}$. Then, the sequence $\{\dbP^N \circ (\m^N(\bX))^{-1} \}_{N \ge 1}$ is tight, and all its accumulation points are supported on $\cP_L(t,\m)$.
\end{prop}



\begin{defn}
{\rm (i)} Denote $Y := (A,M)$ the canonical process on $\O^2 := \O \times \O$. Let $\tilde \cP_L$ be the set of probability measures $\dbP$ on $\O^2$ such that:\\
$\bullet$ $A$ is absolutely continuous w.r.t. to the Lebesgue measure on $[0,T]$, with $\lvert \frac{dA_s}{ds} \rvert \le L$, $\dbP$-a.s.,\\
$\bullet$ $M$ is a $\dbP$-martingale on $[0,T]$, with $\sqrt{\frac{d\langle M \rangle s}{ds}}  \le L$, $\dbP$-a.s.

\no{\rm (ii)} Denote $\bY := (\bA,\bM) = \{(A^k, M^k)\}_{k \in [N]}$ the canonical process on $\O^{N,2} := \O^N \times \O^N$. Let $\tilde \cP_L^N$ be the set of probability measures $\dbP$ on $\O^{N,2}$ s.t.:\\
$\bullet$ $\bA$ is absolutely continuous w.r.t. to the Lebesgue measure on $[0,T]$, with $\lvert \frac{dA^k_s}{ds} \rvert \le L$, for all $k \in [N]$, $\dbP$-a.s.,\\ 
$\bullet$ $\bM$ is a $\dbP$-martingale on $[0,T]$, with $ \sqrt{\frac{d\langle M^k \rangle s}{ds}}  \le L$ and $\langle M^k, M^l \rangle = 0$, for all $k \neq l \in [N]$, $\dbP$-a.s.
\end{defn}

Since a semimartingale is defined as the sum of a finite variation process $\bA$ and a local martingale $\bM$, it is more convenient to show first the tightness of the sequence of empirical measures associated with the pair $(\bA, \bM)$ rather than handling the sum $\bA + \bM$ directly, as it is simpler to show that their properties ``propagate" independently. 

\begin{lem}\label{lem-propagation1}
For all $\{\dbP^N \in \tilde \cP_L^N\}_{N \ge 1}$, the sequence $\{\dbP^N \circ (\m^N(\bY))^{-1} \}_{N \ge 1}$ is tight, and all its accumulation points are supported on $\tilde \cP_L$.
\end{lem}

\proof
{\bf Step 1:} We first prove the existence of a converging subsequence. For all $N \ge 1$, denote $\n^N := \dbP^N \circ (\m^N(\bY))^{-1} \in \cP(\cP_2(\O^2))$. By Lacker \cite[Corollary B.1]{LackerMartingale}, we have to prove that \\
{\rm (i)} $\{ \n^N \}_{N \ge 1}$ is uniformly integrable, i.e., $\lim_{R \to \infty} \sup_{N \ge 1} \dbE^{\n^N}\Big[ \cW_2^2(\l, \d_0)\1_{\cW_2(\l, \d_0) \ge R} \Big] = 0$, where $\l$ is the identity map on $\cP_2(\O^2)$. \\
{\rm (ii)} the sequence of mean measures $\{ \dbE^{\dbP^N}\big[\m^N(\bY) \big]\}_{N \ge 1}$ is tight, \\
where, for all  $\dbP \in \cP_2(\O^{N,2})$ and $\tilde \m : \O^{N,2} \longrightarrow \cP_2(\O^2)$, the mean measure $\dbE^{\dbP}[\tilde \m] \in \cP(\O^2)$ is defined by
$$ \big\langle \dbE^{\dbP}[\tilde \m], \f \big\rangle := \dbE^\dbP \big[ \langle \tilde \m, \f \rangle \big] \q \mbox{for all $\f \in C_b^0(\O^2)$}. $$
Let $R > 0$. We have
\beaa
\dbE^{\n^N}\Big[ \cW_2(\l, \d_0) \1_{\cW_2(m, \d_0)  \ge R}  \Big] &=& \dbE^{\dbP^N}\Big[ \cW_2\big(\m^N(\bY), \d_0\big) \1_{\cW_2(\m^N(\bY), \d_0) \ge R}  \Big] \\
&\le& \frac 1 R  \dbE^{\dbP^N}\Big[ \cW_2^2 \big(\m^N(\bY), \d_0\big)  \Big] \le \frac 2 R  \dbE^{\dbP^N}\Big[ \frac 1 N \sum_{i=1}^N \lvert A^i \rvert^2 + \lvert M^i \rvert^2  \Big]. \\
\eeaa
For each $i \in [N]$, we have 
$ \dbE^{\dbP^N}\big[\lvert A^i \rvert^2 \big] \le (LT)^2 \ \mbox{and} \ \dbE^{\dbP^N}\big[\lvert M^i \rvert^2 \big] \le 4L^2T, $
the latter by Doob's inequality. Thus, there exists a constant $C_{T,L}$ independent from $N$ and $R$ s.t.
$$ \dbE^{\n^N}\Big[ \cW_2(\l, \d_0) \1_{\cW_2(\l, \d_0)  \ge R}  \Big] \le \frac{C_{T,L}}{R} \q \mbox{for all $N \ge 1$ and $R \ge 0$,} $$
and therefore $\lim_{R \to \infty} \sup_{N \ge 1}  \dbE^{\n^N}\Big[ \cW_2(\l, \d_0) \1_{\cW_2(\l, \d_0)  \ge R}  \Big] = 0$ and (i) is proved.

To show that $\{ \dbE^{\dbP^N}\big[\m^N(\bY) \big]\}_{N \ge 1}$ is tight, we prove Aldous' criterion (see Billingsley \cite[Theorem 16.10]{Bil}), i.e.,
\bea\label{aldous}
\sup_{N \ge 1} \sup_{\t \in \cT_{0,T}} \Big\langle \dbE^{\dbP^N}\big[\m^N(\bY) \big], \lvert A_{(\t+\d) \wedge T} - A_\t \rvert^2 + \lvert M_{(\t+\d) \wedge T} - M_\t \rvert^2 \Big\rangle \underset{\d \to 0}{\longrightarrow} 0,
\eea
where $\cT_{0,T}$ denotes the set of $[0,T]$-valued $\dbF$-stopping times. Yet, for fixed $N, \t$ and $\d$,
\beaa
\Big\langle \dbE^{\dbP^N}\big[\m^N(\bY) \big], \lvert M_{(\t+\d) \wedge T} - M_\t \rvert^2 \Big\rangle  =   \frac 1 N \sum_{i=1}^N \dbE^{\dbP^N}\Big[ \big\lvert M_{(\t+\d) \wedge T} ^i - M_\t^i \big\rvert^2 \Big] \le L^2 \d
\eeaa
by Itô's isometry. We obtain a similar estimate for $A$, and this implies \eqref{aldous} and consequently (ii), and thus $\{\n^N\}_{N \ge 1}$ admits a subsequence converging to some $\n \in \cP(\cP_2(\O^2))$.\\
\no{\bf Step 2:} We show that all the accumulations points are supported on $\tilde \cP_L$, i.e.\ that $\n$ is supported on $\tilde \cP_L$. Observe that, by definition of $\dbP^N$, we have
$$ \lvert A^k_s - A^k_r \rvert \le L\lvert s - r \rvert, \ \mbox{$\dbP^N$-a.s., for all $k \in [N]$ and $s,r \in [0,T]$,} $$
and thus 
$$ \lvert A_s - A_r \rvert \le L\lvert s - r \rvert, \ \mbox{$\m^N(\bY)$-a.s., $\dbP^N$-a.s., for all $k \in [N]$ and $s,r \in [0,T]$,}  $$
and finally
$$ \n^N\Big[ \l\Big( \lvert A_s - A_r \rvert \le L\lvert s - r \rvert \Big) = 1 \Big] = 1, \ \mbox{for all $N \ge 1$ and $s,r \in [0,T]$.}  $$
Since $\{ \lvert A_s - A_r \rvert \le L\lvert s - r \rvert \}$ is closed in $\O^2$, $\Big\{ \l\Big( \lvert A_s - A_r \rvert \le L\lvert s - r \rvert \Big) = 1 \Big\}$ is closed in $\cP_2(\O^2)$, and thus the weak convergence of $\n^N$ to $\n$ implies
$$ 1 = \underset{N \to \infty}{\lim \sup} \ \n^N\Big[ \l\Big( \lvert A_s - A_r \rvert \le L\lvert s - r \rvert \Big) = 1 \Big] \le \n\Big[ \l\Big( \lvert A_s - A_r \rvert \le L\lvert s - r \rvert \Big) = 1 \Big] \le 1, $$
that is, $\n\Big[ \l\Big( \lvert A_s - A_r \rvert \le L\lvert s - r \rvert \Big) = 1 \Big] = 1$. Since $s$ and $r$ are arbitrary, this implies that $A$ is absolutely continuous w.r.t. the Lebesgue measure on $[0,T]$ with $\lvert \frac{dA_s}{ds} \rvert \le L$, $\l$-a.s., $\n$-a.s. 

We now prove that $M$ is a $\l$-martingale on $[0,T]$, $\n$-a.s. Fix $r \le s$ in $[0,T]$, and $h_r := h(Y_r)$, where $h \in C_b^0(\O^2)$. We compute:
\beaa
\dbE^{\n^N}\Big[ \big\langle \l, h_r(M_s - M_r) \big\rangle^2 \Big] &=& \dbE^{\dbP^N}\Big[ \Big( \frac 1 N \sum_{i=1}^N h(Y_r^i)(M_s^i - M_r^i)\Big)^2 \Big] \\
&\le& \frac{1}{N^2} \sum_{i=1}^N  \lvert h \rvert^2 \dbE^{\dbP^N}\Big[ \lvert M_s^i - M_r^i \rvert^2 \Big] 
\le \frac{\lvert h \rvert^2 L^2 T}{N} \underset{N \to \infty}{\longrightarrow} 0,
\eeaa
where we used the fact that $\langle M^k, M^l \rangle \1_{k \neq l} = 0$, and the $\si(Y_r)$-measurability of $h_r$ to derive the first inequality. Thus, as $\n^N$ converges weakly to $\n$,
$$ 0 \le  \dbE^{\n}\Big[ \big\langle \l, h_r(M_s - M_r) \big\rangle^2 \Big] \le \underset{N \to \infty}{\lim \inf} \ \dbE^{\n^N}\Big[ \big\langle \l, h_r(M_s - M_r) \big\rangle^2 \Big] = 0,$$
hence $\dbE^{\n}\Big[ \big\langle \l, h_r(M_s - M_r) \big\rangle^2 \Big] = 0$, which implies that
$$ \big\langle \l, h_r(M_s - M_r) \big\rangle = 0, \ \mbox{$\n$-a.s.,} $$
which by the arbitrariness of $s,r$ and $h$ means that $M$ is a $\l$-martingale, $\n$-a.s. We prove similarly to $A$ that $\sqrt{\frac{d\langle M \rangle_s}{ds}} \le L$, $\l$-a.s., $\n$-a.s.
\qed

We eventually prove Proposition \ref{lem-propagation} by deriving the tightness of the processes $\{\m^N(\bA + \bM) \}_{N \ge 1}$ from the one of the processes $\{ \m^N(\bA, \bM) \}_{N \ge 1}$.

\no{\bf{Proof of Propostion \ref{lem-propagation}}}
Introduce 
$$ \bA_s^N := \bX_t + \1_{s \ge t}\int_t^s b_r^{\dbP^N}dr, \q \bM_s^N := \1_{s \ge t}\int_t^s \si_r^{\dbP^N}dW_r^{\dbP^N}, $$
where $b^{\dbP^N}, \si^{\dbP^N}$ and $W^{\dbP^N}$ are as in \eqref{cPLN}. Then, we clearly have
$$ \tilde \dbP^N := \dbP^N_{(\bA^N, \bM^N)} \in \tilde \cP_L^N.$$
Therefore, by Lemma \ref{lem-propagation1}, $\n^N := \tilde \dbP^N \circ (\m^N(\bY))^{-1}$ converges weakly to some $\n$ supported on $\tilde \cP_L(t,\m)$. Define $\hat \m^N := \dbP^N \circ (\m^N(\bX))^{-1}$ and fix $\f \in C_b^0(\cP_2(\O))$. We have:
\beaa
\langle \hat \m^N, \f \rangle &=& \dbE^{\dbP^N}\Big[\f\big(\m^N(\bX)\big)\Big] = \dbE^{\dbP^N}\Big[\f\big(\m^N(\bY) \circ (A+M)^{-1}\big)\Big] \\
&=& \dbE^{\n^N}\Big[\f\big(\l \circ (A+M)^{-1}\big)\Big] \underset{N \to \infty}{\longrightarrow} \dbE^{\n}\Big[\f\big(\l \circ (A+M)^{-1}\big)\Big] 
\eeaa
by weak convergence of $\n^N$ to $\n$, since $\l \mapsto \f\big(\l \circ (A+M)^{-1}\big) \in C_b^0(\cP_2(\O^2))$. Thus we have
$$ \langle \hat \m^N, \f \rangle \underset{N \to \infty}{\longrightarrow} \langle \hat \m, \f \rangle, $$
where $\hat \m := \n \circ \big( \l \circ (A+M)^{-1} \big)^{-1}$. Therefore, by arbitrariness of $\f$, $\m^N$ converges weakly to $\hat \m$, which is clearly supported on $\cP_L(t,\m)$ as $\n$ is supported on $\tilde \cP_L$ and $\m^N(\bo^N) \overset{\cW_2}{\longrightarrow} \m$.
\qed

\appendix

\section{Technical lemma}

\begin{lem}\label{lem:convex}
Fix $\bx \in \dbR^{d \times N}$, and introduce for all $h > 0$:
$$ D_h := \{ \by \in \dbR^{d \times N} : \cW_2(\m^N(\by), \m^N(\bx)) \le h \}. $$
Then, for $h$ sufficiently small, $D_h$ is a disjoint union of convex sets. 
\end{lem}
\proof
By Birkhoff’s theorem (see e.g.\ Villani \cite[p.\ 5]{Villani}), we have for all $\by$:
$$ \cW_2(\m^N(\by), \m^N(\bx)) = \min_{\pi \in \mathfrak{S}_N} \ \lVert \by - \pi(\bx) \rVert_2, $$
where $\mathfrak{S}_N$ is the symmetric group of $[N]$, and where we denote by $\pi(\bx) := (x_{\pi(1)}, \dots, x_{\pi(N)})$. Then we have:
$$ D_h = \{ \by \in \dbR^{d \times N} : \exists \pi \in \mathfrak{S}_N, \ \mbox{s.t.} \ \cW_2(\m^N(\by), \m^N(\pi(\bx))) \le h \}, $$
and $D_h$ writes therefore as a finite union of balls with centers in $\{\pi(\bx)\}_{\pi \mathfrak{S}_N}$. Thus, for $h$ small enough, $D_h$ is a disjoint union of balls, which are convex sets. 
\qed

{\small
\bibliographystyle{plain}
\bibliography{references2}}
\end{document}